\documentclass[preprint,12pt]{elsarticle}

\usepackage{amsthm}
\newtheorem{theorem}{Theorem}[section]
\newtheorem{proposition}[theorem]{Proposition}
\newtheorem{remark}{Remark}[section]
\newtheorem{definition}[theorem]{Definition}

\usepackage{amsmath}
\usepackage{amssymb}	
\usepackage{bbold}
\usepackage{bm} 

\usepackage{xcolor}
\definecolor{myBlue}{RGB}{50,117,180}
\definecolor{myRed}{RGB}{200,40,40}
\definecolor{myGreen}{RGB}{34,139,34}
\definecolor{myLightBlue}{RGB}{193,223,255}
\definecolor{myOrange}{RGB}{255,185,0}

\usepackage{cancel} 
\usepackage[normalem]{ulem} 
\usepackage{changepage}
\usepackage{float}
\usepackage{mathrsfs}
\usepackage{url}
\usepackage{xfrac}
\usepackage{graphicx}
\usepackage{multirow}
\usepackage{lipsum} 
\graphicspath{{/}}

\usepackage[font=small,labelfont=bf,format=plain,singlelinecheck=false]{caption} 

\newcommand{\comment}[1]{\ignorespaces} 

\newcommand{\w}[1]{\ensuremath{\mathbf{#1}}} 
\newcommand{\gv}[1]{\ensuremath{\mbox{\boldmath$ #1 $}}} 
\newcommand{\pd}[2]{\frac{\partial #1}{\partial #2}} 

\DeclareMathAlphabet{\mathsfit}{\encodingdefault}{\sfdefault}{m}{n}
\SetMathAlphabet{\mathsfit}{bold}{\encodingdefault}{\sfdefault}{bx}{n}
\newcommand{\abs}[1]{\left| #1 \right|} 
\newcommand{\ws}[1]{\mathsf{#1}} 

\makeatletter
\newcommand*\bigcdot{\mathpalette\bigcdot@{.5}}
\newcommand*\bigcdot@[2]{\mathbin{\vcenter{\hbox{\scalebox{#2}{$\m@th#1\bullet$}}}}}
\makeatother

\makeatletter
\newcommand{\raisemath}[1]{\mathpalette{\raisem@th{#1}}}
\newcommand{\raisem@th}[3]{\raisebox{#1}{$#2#3$}}
\makeatother

\usepackage{accents}

\usepackage{stackengine} 

\makeatletter 
\renewcommand*\env@matrix[1][*\c@MaxMatrixCols c]{%
  \hskip -\arraycolsep
  \let\@ifnextchar\new@ifnextchar
  \array{#1}}
\makeatother

\usepackage{tikz}
\usetikzlibrary{cd}
\makeatletter
\providecommand\pgf@matrix@last@nextcell@options{}
\makeatother

\usepackage{pgfplots}
\usepgfplotslibrary{patchplots}
\usepgfplotslibrary{external}
\usepackage{tikz}
\usetikzlibrary{external}
\usetikzlibrary{calc}
\usetikzlibrary{arrows,shapes,backgrounds,arrows.meta,shapes.arrows,automata,quotes,intersections}
\usetikzlibrary{decorations.text}
\usetikzlibrary{decorations.fractals}
\usetikzlibrary{decorations.pathmorphing}
\usetikzlibrary{decorations.markings}

\usetikzlibrary{decorations.pathreplacing}
\makeatletter
    \let\pgf@decorate@@brace@brace@code@old\pgf@decorate@@brace@brace@code
    \def\pgf@decorate@@brace@brace@code{
        \ifdim\pgfdecoratedremainingdistance<4\pgfdecorationsegmentamplitude
            \pgftransformxscale{\pgfdecoratedremainingdistance/4\pgfdecorationsegmentamplitude}
            \pgfdecoratedremainingdistance=4\pgfdecorationsegmentamplitude
        \fi
        \pgf@decorate@@brace@brace@code@old
    }
\makeatother

\newcommand{\norm}[1]{\left\lVert#1\right\rVert}

\newcommand*\mA{\mathbb{A}}

\newcommand*\mC{\mathbb{C}}

\newcommand*\mmD{\mathbb{D}}

\newcommand*\mG{\mathbb{G}}

\newcommand*\mI{\mathcal{I}}
\newcommand*\mJ{\mathcal{J}}
\newcommand*\mK{\mathbb{K}}

\newcommand*\mM{\mathbb{M}}

\newcommand*\mN{\mathbb{N}}
\newcommand*\mO{\mathcal{O}}
\newcommand*\mQ{\mathbb{Q}}
\newcommand*\mP{\mathbb{P}}
\newcommand*\mmP{\mathcal{P}}
\newcommand*\mPi{\mathbb{\Pi}}
\newcommand*\mR{\mathbb{R}}
\newcommand*\mS{\mathbb{S}}

\newcommand*\mT{\mathcal{T}}
\newcommand*\mmT{\mathbb{T}}

\newcommand*\mU{\mathbb{U}}
\newcommand*\mV{\mathbb{V}}
\newcommand*\mW{\mathcal{W}}
\newcommand*\mmW{\mathbb{W}}

\newcommand*\mOne{\mathbb{1}}
\newcommand*\mZero{\mathbb{0}}

\newcommand*\md{\mathrm{d}}

\newcommand*\mmOmega{\Omega}
\DeclareMathAlphabet{\mathpzc}{OT1}{pzc}{m}{it}

\usepackage[most]{tcolorbox}

\makeatletter
\newcommand*{\defeq}{\mathrel{\rlap{%
                     \raisebox{0.3ex}{$\m@th\cdot$}}%
                     \raisebox{-0.3ex}{$\m@th\cdot$}}%
                     =}
\newcommand*{\eqdef}{=\mathrel{\rlap{%
                     \raisebox{0.3ex}{$\m@th\cdot$}}%
                     \raisebox{-0.3ex}{$\m@th\cdot$}}%
                     }
\makeatother

\usepackage{mathtools}

\usepackage{hhline}

\usepackage{stmaryrd}

\usepackage{enumitem}

\definecolor{light-gray}{gray}{.85}
\newsavebox{\songboxbox}

\usepackage[most]{tcolorbox}
\newtcbox{\MyBox}[1][red]{on line, size=tight, boxsep=1pt, colframe=#1!50!black, colback=#1!10!white}

\usepackage[framemethod=TikZ]{mdframed}
\mdfsetup{%
    outerlinewidth=1pt,
    innertopmargin=6pt,
    innerbottommargin=6pt,
    skipabove=10pt,
    skipbelow=10pt,
    backgroundcolor=black!10,
    roundcorner=20pt
}

\DeclareFontFamily{U}{matha}{\hyphenchar\font45}
\DeclareFontShape{U}{matha}{m}{n}{
      <5> <6> <7> <8> <9> <10> gen * matha
      <10.95> matha10 <12> <14.4> <17.28> <20.74> <24.88> matha12
      }{}
\makeatletter
\newcommand{\blandor}[1]{\mathbin{\@blandor{#1}}}
\newcommand{\@blandor}[1]{\mathchoice
  {\@@blandor{#1}{\tf@size}}
  {\@@blandor{#1}{\tf@size}}
  {\@@blandor{#1}{\sf@size}}
  {\@@blandor{#1}{\ssf@size}}
}
\newcommand{\@@blandor}[2]{%
    \raisebox{.1ex}{\rotatebox[origin=c]{#1}{%
      \fontsize{#2}{#2}\usefont{U}{matha}{m}{n}\symbol{\string"CE}}}%
}
\makeatother

\usepackage{slashed} 
\newcommand{\cmmnt}[1]{\ignorespaces} 

\usepackage{environ}
\NewEnviron{eqn}{
\begin{align}\begin{split}
\BODY
\end{split}\end{align}
}

\usepackage{xfrac} 

\usepackage{bbding}
\DeclareFontFamily{U}{MnSymbolC}{}
\DeclareSymbolFont{MnSyC}{U}{MnSymbolC}{m}{n}
\DeclareMathSymbol{\diamondplus}{\mathbin}{MnSyC}{"7C}
\DeclareMathSymbol{\diamonddot}{\mathbin}{MnSyC}{"7E}
\DeclareFontShape{U}{MnSymbolC}{m}{n}{
    <-6>  MnSymbolC5
   <6-7>  MnSymbolC6
   <7-8>  MnSymbolC7
   <8-9>  MnSymbolC8
   <9-10> MnSymbolC9
  <10-12> MnSymbolC10
  <12->   MnSymbolC12}{}
  
  \newcommand\underlay[4]{%
  \stackengine{0pt}%
  {\kern#2\includegraphics[height=#1]{#4}}%
  {\includegraphics[height=#1]{#3}}%
  {O}{l}{F}{F}{L}%
}
\newcommand\addunderlay[4]{%
  \stackengine{0pt}%
  {\kern#2\includegraphics[height=#1]{#4}}%
  {#3}%
  {O}{l}{F}{F}{L}%
}

\usepackage{hyperref}

\begin{document}

\begin{frontmatter}

\title{Generalized Yee methods:\\Scalable symplectic finite element Maxwell solvers}


\author[1]{Alexander S. Glasser\corref{cor1}}
\ead{asg5@princeton.edu}
\author[1]{Hong Qin}
\ead{hongqin@princeton.edu}

\affiliation[1]{Department of Astrophysical Sciences, Princeton University, Princeton, NJ 08544, USA}

\cortext[cor1]{Corresponding author}


\begin{abstract}
Yee's finite-difference method preserves two crucial properties of Maxwell's equations---locality and symplecticity---and thereby enjoys two computational advantages: scalability on high-performance architectures and long-time numerical accuracy. In this work, we show that Yee's method is a special case of a class of structure-preserving finite element methods---termed \textit{generalized Yee methods} (GYMs)---that are designed to retain both crucial properties. GYMs are built from de Rham-conforming finite elements and achieve locality through sparse mass matrices and their sparse approximate inverses (SPAIs). We prove that the symplectic structure of GYMs is invariant under such sparse approximations, freeing the choice of sparsification strategy. We introduce a novel sparsification strategy, SPAI-OP, which concentrates accuracy at prescribed wave modes by \emph{operator probing}. We further extend GYMs to structure-preserving electromagnetic particle-in-cell (PIC) methods, whose symplecticity over particle trajectories requires the smooth fields afforded by higher-order finite elements. GYMs therefore retain the computational virtues of Yee's method while enabling unstructured meshes, higher-order accuracy, spectral adaptivity, and symplectic particle coupling.
\end{abstract}

\begin{keyword}
Generalized Yee methods \sep
structure-preserving algorithms \sep
finite element exterior calculus \sep
operator-optimized sparse approximate inversion \sep
symplectic integrator \sep
electromagnetic simulation
\end{keyword}



\end{frontmatter}

\section{Introduction}

The Yee algorithm \cite{yee_numerical_1966,taflove_computational_2005}---alternatively, the finite difference time domain (FDTD) method---defines electromagnetic fields on a cubical mesh. It associates to each \emph{edge} a component of the electric field $\w{E}$, and to each \emph{face} a component of the magnetic field $\w{B}$. This discretization reflects a natural geometric description of Maxwell's equations, in which ${\w{E}\in\Lambda^1(\mR^3)}$ is a differential $\text{1-form}$ and ${\w{B}\in\Lambda^2(\mR^3)}$ is a differential $\text{2-form}$. In this way, Yee's method employs a technique adopted in many \emph{structure-preserving algorithms} \cite{hairer_geometric_2006}, wherein differential $k$-forms are discretized by associating them with $k$-dimensional features of a mesh \cite{whitney_geometric_1957,desbrun_discrete_2005,arnold_finite_2006,arnold_finite_2010}.

Structure-preserving algorithms have been widely adopted in computational physics, including in gravitational simulations  \cite{gladman_symplectic_1991,kinoshita_symplectic_1991,chambers_symplectic_2002,bravetti_numerical_2020,kur_discrete_2022}, geophysics \cite{li_structure-preserving_2012,liu_modified_2015} and plasma physics \cite{squire_geometric_2012,xiao_explicit_2015,he_hamiltonian_2015,crouseilles2015Hamiltonian,qin_canonical_2016,kraus_gempic:_2017,morrison_structure_2017,glasser_gauge-compatible_2022}. Such algorithms generally derive from variational principles or Hamiltonian systems, and as a result preserve essential mathematical features of their underlying physical systems, including symplectic structure, topology, symmetries, and conservation laws. Despite Yee's omission of any overt Lagrangian or Hamiltonian formulations \cite{yee_numerical_1966}, Yee's method is (apparently serendipitously) one of the most historically successful examples of a structure-preserving algorithm \cite{stern_geometric_2009}.

Restricting our discussion to the symplectic Hamiltonian setting, the essential appeal of structure-preserving algorithms is that their trajectories over a fixed step size integrate a Hamiltonian `near' the true Hamiltonian $H$. Specifically, an order $p$ algorithm evolves exactly along the trajectory of a Hamiltonian $H_\epsilon$ that is $\epsilon^p$-close to $H$ \cite{benettin_hamiltonian_1994}, i.e.\ ${\norm{H-H_\epsilon}\leq C\epsilon^p}$ for a constant $C$. In this sense, the simulation's degrees of freedom (DOFs) traverse a physically realizable manifold near the true one. This accuracy is generally more valuable than energy conservation (${\dot{H}=0}$) alone: whereas energy-conserving algorithms place a single scalar constraint on their $N$ DOFs, symplectic algorithms satisfy $N$ constraints by preserving the entire hierarchy of Poincar\'e integral invariants \cite{kraus_gempic:_2017}. Furthermore, while phase errors may still secularly accumulate in symplectic integrators \cite{scuro_forward_2005}, they do so at slower rates than in many non-symplectic integration algorithms \cite{gladman_symplectic_1991}. Overall, symplectic algorithms can offer a significant advantage in simulation accuracy for a wide range of physical phenomena.

Well-chosen finite element families are readily purposed in the construction of symplectic algorithms. To make a finite element method structure-preserving, one requires that the projection maps ${\Pi^k:H\Lambda^k\rightarrow\Lambda^k(\mT_h)}$---from Sobolev spaces $H\Lambda^k$ to their finite element subspaces $\Lambda^k(\mT_h)$---commute with the exterior derivative $\md$, that is: ${\md\circ\Pi^i=\Pi^{i+1}\circ\md}$. (See the commutative diagram of Fig.~\ref{deRhamDiagram}.) The development of finite element exterior calculus (FEEC) \cite{arnold_finite_2006,arnold_finite_2010} firmly established the critical importance of this commuting property in ensuring stable and well-posed finite element methods.

\begin{figure}[b!]
\begin{tikzcd}[column sep=0.45in, row sep=0.8in]
0 \arrow[r] & H\Lambda^0 \arrow[r,"\md"] \arrow[d,"\Pi^0"] & H\Lambda^1 \arrow[r,"\md"] \arrow[d,"\Pi^1"] & H\Lambda^2 \arrow[r,"\md"] \arrow[d,"\Pi^2"] & H\Lambda^3 \arrow[r] \arrow[d,"\Pi^3"] & 0 \\
0 \arrow[r] & \Lambda^0(\mT_h) \arrow[r,"\md"] & \Lambda^1(\mT_h) \arrow[r,"\md"] & \Lambda^2(\mT_h) \arrow[r,"\md"] & \Lambda^3(\mT_h) \arrow[r] & 0
\end{tikzcd}
\caption{The commuting diagram for derivative and projection mappings between Sobolev spaces ${H\Lambda^k}$ and their de~Rham-conforming finite element subspaces ${\Lambda^k(\mT_h)}$ on a mesh ${\mT_h\subset\mR^3}$ of maximal cell diameter $h$.}
\label{deRhamDiagram}
\end{figure}

This commutation guarantees the topological fidelity of a finite element algorithm. To see this, note that the sequence in the bottom row of Fig.~\ref{deRhamDiagram} is a chain complex of finite element subspaces, since ${\Lambda^k(\mT_h)\subset H\Lambda^k}$ and ${\md\circ\md=0}$ between Sobolev spaces. The commuting projections $\Pi^p$ ensure that this discrete finite element chain complex inherits the cohomology of the continuous (Sobolev) chain complex in the top row. This commutation has also been shown critical in establishing gauge symmetry and conservation laws, as demonstrated by the conservation of charge via Gauss' law in plasma physics algorithms \cite{squire_geometric_2012,xiao_explicit_2015,kraus_gempic:_2017,bossavit_whitney_1988,glasser_geometric_2020,glasser_gauge_2022}.

A judicious treatment of finite element methods also readily assimilates the other crucial feature of Yee's method---scalability. Scalable finite element Yee-type methods have proliferated in the literature over the last few decades, using a wide array of approaches to generate sparse approximations of mass matrices and their inverses, including: mass lumping \cite{cohen_gauss_1998,egger_mass_lumped_2020}, thresholding off-diagonal elements \cite{bo_he_sparse_2006, he_differential_2007}, optimizing a fixed-sparsity-pattern inverse mass matrix \cite{kim_parallel_2011,teixeira_differential_2013}, and broken finite elements that overlap only cell-wise \cite{campos_pinto_gauss-compatible_2016,kapidani_arbitrary-order_2021,guclu_broken_2023}. As we will show, all of these sparsification methods are compatible with preserving the symplecticity of GYMs.

In this paper, we organize and extend these methods into a single class, which we call \emph{generalized Yee methods} (GYMs). A GYM is a time-domain method, whose time discretization is defined by a Hamiltonian splitting method, and whose spatial discretization is defined by finite elements whose mass matrices and their (approximate) inverses are required to be sparse. The finite element families employed are further required to be de~Rham-conforming, in the sense of commutation in Fig.~\ref{deRhamDiagram}. We show that Yee's method is recovered as the simplest possible GYM, in which a Strang splitting \cite{strang_construction_1968} evolves electromagnetic fields on a cubical mesh, whose finite element spaces are Whitney forms with mass matrices approximated as diagonal.

The central theoretical contribution of this work is a proof that the symplecticity of GYMs is independent of the choice of mass-matrix approximation: any approximate mass-matrix inverse substituted into the finite element Hamiltonian yields a symplectic algorithm. We further require that approximate mass matrices and their inverses be symmetric positive-definite (SPD)---not to preserve symplecticity, but to ensure numerical stability. This result unifies the literature on scalable Yee-type methods within a single structure-preserving framework, and, crucially, shows that sparsification may be designed with considerable flexibility: once symmetry and positive-definiteness are enforced, accuracy and computational efficiency become the primary design criteria.

This freedom motivates a second contribution of the paper: an operator-probed sparse approximate inverse (SPAI-OP) method that concentrates accuracy on prescribed wave modes of a GYM's discrete curl-of-curl operator. Inspired by the modified SPAI (MSPAI) framework of Huckle and Kallischko \cite{huckle_frobenius_2007}, who introduced probing constraints to improve sparse approximate inverses on targeted subspaces, SPAI-OP augments the standard SPAI objective with soft penalties that encourage the approximate inverse to be accurate on user-specified wave modes (e.g. operator eigenmodes, specific initial conditions, or wave packets at frequencies of physical interest). The proposed method thereby serves as a concrete example of the broader sparsification latitude enabled by the symplecticity result. The symmetry-constrained optimality condition takes the form of a generalized Sylvester system, solvable by a scalable matrix-free preconditioned conjugate gradient method. The approach is distinct from prior spectral optimizations---such as dispersion-optimized finite-difference stencils \cite{cole_high_1997,cowan_generalized_2013,blinne_systematic_2018}---in that it optimizes the inverse mass matrix \cite{grote_spai_1997} within a general finite element framework, preserving the exact topological curl operator inherited from the de~Rham complex (Fig.~\ref{deRhamDiagram}).

These results are particularly consequential for kinetic plasma simulations.  Structure-preserving electromagnetic PIC methods have been an active area of development within the FEEC framework \cite{squire_geometric_2012,kraus_gempic:_2017,glasser_gauge-compatible_2022}, but have typically relied upon dense mass-matrix inversion, limiting their scalability. Recently, Barham and Burby \cite{barham_diagnosing_2025} demonstrated that symplecticity over particle trajectories requires discrete electromagnetic fields that are pointwise $C^1$---a condition violated by most finite element families and by standard finite differences. Among de~Rham-conforming finite element families, pointwise smoothness is provided by B-splines of degree ${r\geq2}$ on cubical meshes \cite{kraus_gempic:_2017,barham_diagnosing_2025,buffa_isogeometric_2010,buffa_isogeometric_2011}. For simplicial meshes, smooth de~Rham complexes have been constructed using Clough-Tocher macro-elements and related techniques \cite{neilan_smooth_2015,christiansen_nodal_2018,christiansen_generalized_2018,fu_alfeld_2020}. As a result, GYMs are extensible to particle-coupled PIC codes as well.

The combination of these results establishes GYMs as a practical framework that retains the defining computational virtues of Yee's method---scalability and symplecticity---while admitting unstructured meshes, higher-order accuracy, spectral adaptivity, and symplectic particle coupling.

The remainder of this article is organized as follows. Section~\ref{FEStructGYMs} reviews the FEEC formalism and its application to Maxwell's equations. Section~\ref{SymplecticStructGYMs} defines the GYM class and proves that its symplecticity is independent of mass-matrix approximation (Theorem~\ref{thm:symplecticity}), while stability requires only positive definiteness (Proposition~\ref{prop:stability}). Section~\ref{YeeIsFEEC} recovers Yee's algorithm as the simplest GYM. Section~\ref{SparseMassMatSect} reviews existing SPAI techniques and introduces the SPAI-OP formulation. Section~\ref{PICsect} extends GYMs to structure-preserving PIC methods, addressing the Barham--Burby $C^1$ smoothness requirement. Numerical results are presented in Section~\ref{NumResult}, and Section~\ref{ConclusionSect} concludes.

\section{Finite Element Structure of GYMs\label{FEStructGYMs}}

We briefly establish notation from finite element exterior calculus (FEEC) \cite{arnold_finite_2006,arnold_finite_2010}; see \cite{kraus_gempic:_2017,glasser_gauge-compatible_2022} for its application to electromagnetic simulations. Let $\Lambda^p(\mT_h)$ denote the space of finite element differential $p$-forms on a simplicial or cubical complex ${\mT_h\subset\mR^n}$ of maximal cell diameter $h$, spanned by a locally supported finite element basis ${\gv{\Lambda}^p}$ of dimension $N_p$. An arbitrary $p$-form ${\w{S}\in\Lambda^p(\mT_h)}$ is written ${\w{S}(\w{x})=\w{s}\cdot\gv{\Lambda}^p(\w{x})=s_i\Lambda^p_i(\w{x})}$ with ${\w{s}\in\mR^{N_p}}$ (Einstein summation hereafter). Individual spatial components are denoted with Greek letters for coordinate indices, as in ${\w{S}(\w{x})_{\mu_1\cdots\mu_p}=s_i\Lambda^p_i(\w{x})_{\mu_1\cdots\mu_p}}$.

Since $\md$ is linear, it is represented by matrices---the gradient ${\mG\in\mR^{N_1\times N_0}}$, curl ${\mC\in\mR^{N_2\times N_1}}$, and divergence ${\mmD\in\mR^{N_3\times N_2}}$---defined such that ${\mG^T\gv{\Lambda}^1=\md\gv{\Lambda}^0}$, ${\mC^T\gv{\Lambda}^2=\md\gv{\Lambda}^1}$, ${\mmD^T\gv{\Lambda}^3=\md\gv{\Lambda}^2}$, and satisfying ${\mC\mG=\mZero}$ and ${\mmD\mC=\mZero}$. These comprise the \emph{discrete de~Rham complex},
\begin{eqn}
0\rightarrow\Lambda^0(\mT_h)\xrightarrow{\md~(\cong~\mG)}\Lambda^1(\mT_h)\xrightarrow{\md~(\cong~\mC)}\Lambda^2(\mT_h)\xrightarrow{\md~(\cong~\mmD)}\Lambda^3(\mT_h)\rightarrow 0.
\label{discreteDeRham}
\end{eqn}
A finite element family is called \emph{de~Rham-conforming} if the projections ${\Pi^p:H\Lambda^p\rightarrow\Lambda^p(\mT_h)}$ of Fig.~\ref{deRhamDiagram} commute with $\md$:
\begin{eqn}
\md\circ\Pi^p=\Pi^{p+1}\circ\md,
\label{eqnCommute}
\end{eqn}
ensuring that the discrete complex inherits the cohomology of the continuous one \cite{arnold_finite_2006}. The top row of Fig.~\ref{deRhamDiagram} is the continuous (Sobolev) de~Rham complex, while the bottom row is the discrete finite element complex of Eq.~(\ref{discreteDeRham}). The vertical arrows are the projection maps of Eq.~(\ref{eqnCommute}). Well-known examples of de~Rham-conforming families suitable for GYMs include Whitney forms \cite{whitney_geometric_1957,bossavit_whitney_1988} and N\'ed\'elec finite elements \cite{nedelec_mixed_1980} on simplicial meshes, and B-spline discretizations \cite{buffa_isogeometric_2010,buffa_isogeometric_2011} on cubical meshes.


The mass matrix ${(\mM_p)_{ij}=\int\md\w{x}\,(\Lambda^p_i,\Lambda^p_j)_p}$ encodes the metric structure, with ${(\cdot,\cdot)_p}$ the pointwise inner product on $p$-forms. Each $\mM_p$ is sparse (due to local basis support), symmetric, and positive definite, but its inverse $\mM_p^{-1}$ is typically dense. \emph{A key observation underlying GYMs is that $\mG$, $\mC$, and $\mmD$ encode topology, while the mass matrices encode metric structure.} The mass matrices can therefore be freely approximated without disturbing the topological or---as we prove in Section~\ref{SymplecticStructGYMs}---the symplectic structure.

For electromagnetic simulations, we discretize the magnetic vector potential as ${\w{A}=\w{a}\cdot\gv{\Lambda}^1}$ in the temporal gauge (${\phi=0}$). The magnetic field is then given by
\begin{eqn}
\w{b}\cdot\gv{\Lambda}^2=\w{B}=\md\w{A}=\w{a}\cdot\md\gv{\Lambda}^1=\w{a}\cdot\mC^T\gv{\Lambda}^2=\mC\w{a}\cdot\gv{\Lambda}^2
\label{BfromA}
\end{eqn}
with ${\mC\w{a}=\w{b}\in\mR^{N_2}}$. We define the electric field as ${\w{E}=\w{e}\cdot\mM_1^{-1}\cdot\gv{\Lambda}^1}$, with ${\w{e}\in\mR^{N_1}}$. Following \cite{glasser_gauge-compatible_2022}, we use a convention wherein the coefficients $\w{e}$ determine $\w{E}$ with an additional factor of the exact matrix inverse, $\mM_1^{-1}$. This convention ensures that $(\w{a},\w{e})$ are canonically conjugate, yielding the canonical Poisson bracket of Eq.~(\ref{discHamilSystem}). The pair ${(\w{a},\w{e})\in\mR^{2N_1}}$ defines the DOFs of the finite element electromagnetic fields.

The de~Rham identities guarantee ${\mmD\w{b}=\mmD\mC\w{a}=\w{0}}$ (divergence-free $\w{B}$) and ${\mC\mG=\mZero}$, the latter of which will help endow the discretization with exact charge conservation, as we describe in Section~\ref{PICsect}. These properties, too, follow from the topological structure of the discrete de~Rham complex, and are immune to any approximation of the metric structure via mass matrices.

\section{Symplectic Structure of GYMs\label{SymplecticStructGYMs}}

To describe the dynamics of the discrete electromagnetic fields defined in Section~\ref{FEStructGYMs}, we first recall the canonical Hamiltonian formulation of Maxwell's equations in the continuum, with Poisson bracket and Hamiltonian in SI units \cite{marsden_hamiltonian_1982}:
\begin{eqn}
\{F,G\}_{\text{EM}}&=\frac{1}{\epsilon_0}\int\md\w{x}\left(\frac{\delta F}{\delta \w{E}}\cdot\frac{\delta G}{\delta \w{A}}-\frac{\delta G}{\delta \w{E}}\cdot\frac{\delta F}{\delta \w{A}}\right)\\
H_{\text{EM}}&=\frac{1}{2}\int\md\w{x}\left(\epsilon_0\abs{\w{E}}^2+\frac{1}{\mu_0}\abs{\nabla\times\w{A}}^2\right).
\label{continuumEMHamil}
\end{eqn}
Substituting ${\w{A}=\w{a}\cdot\gv{\Lambda}^1}$ and ${\w{E}=\w{e}\cdot\mM_1^{-1}\cdot\gv{\Lambda}^1}$ yields \cite{glasser_gauge-compatible_2022} the discrete Poisson bracket and Hamiltonian
\begin{eqn}
\{F,G\}_\text{d}&=\frac{1}{\epsilon_0}\left(\pd{F}{\w{e}}\cdot\pd{G}{\w{a}}-\pd{G}{\w{e}}\cdot\pd{F}{\w{a}}\right)\\
H_\text{d}&=\frac{\epsilon_0}{2}\w{e}^T\mM_1^{-1}\w{e}+\frac{1}{2\mu_0}\w{a}^T\mC^T\mM_2\mC\w{a},
\label{discHamilSystem}
\end{eqn}
with equations of motion
\begin{eqn}
\dot{\w{a}}&=\{\w{a},H_\text{d}\}_\text{d}=-\mM_1^{-1}\w{e}\\
\dot{\w{e}}&=\{\w{e},H_\text{d}\}_\text{d}=c^2\mC^T\mM_2\mC\w{a}.
\label{exactEOM}
\end{eqn}
Due to the dense inverse $\mM_1^{-1}$, each timestep of this system generally requires global communication of all degrees of freedom.

The key observation is that the Poisson bracket $\{\cdot,\cdot\}_\text{d}$ is canonical in $(\w{a},\w{e})$---i.e. $(\w{a},\w{e})$ form Darboux coordinates---such that $\{\cdot,\cdot\}_\text{d}$ involves neither $\mM_1$ nor $\mM_2$. The mass matrices enter only through the Hamiltonian. We may therefore respectively replace $\mM_1^{-1}$ (and optionally $\mM_2$) with arbitrary sparse approximations $\mQ_1$ (and optionally $\tilde{\mM}_2$) in the Hamiltonian,
\begin{eqn}
\tilde{H}_\text{d}=\frac{\epsilon_0}{2}\w{e}^T\mQ_1\w{e}+\frac{1}{2\mu_0}\w{a}^T\mC^T\tilde{\mM}_2\mC\w{a},
\label{discTildeHamilSystem}
\end{eqn}
without modifying the symplectic structure. The resulting equations of motion,
\begin{eqn}
\dot{\w{a}}&=-\mQ_1\w{e}\\
\dot{\w{e}}&=c^2\mC^T\tilde{\mM}_2\mC\w{a},
\label{discEOM}
\end{eqn}
are scalable discrete forms of ${\dot{\w{A}}=-\w{E}}$ and ${\dot{\w{E}}=c^2\nabla\times\nabla\times\w{A}}$. (In Eq.~(\ref{discEOM}), $\mQ_1$ and $\tilde{\mM}_2$ are assumed symmetric without loss of generality (WLOG), since only their symmetric parts enter the EOM.) Written equivalently (via Poisson reduction \cite{glasser_geometric_2020}) in terms of ${\w{b}=\mC\w{a}}$, these become
\begin{eqn}
\dot{\w{b}}&=-\mC\mQ_1\w{e}\\
\dot{\w{e}}&=c^2\mC^T\tilde{\mM}_2\w{b},
\label{preYeeFEECEqns}
\end{eqn}
which is the form we shall use to identify Yee's method as a GYM in Section~\ref{YeeIsFEEC}.

We emphasize that the definition ${\w{E}=\w{e}\cdot\mM_1^{-1}\cdot\gv{\Lambda}^1}$, with its dense inverse, remains unmodified by the approximations $\mQ_1$ and ${\tilde{\mM}_2}$; only the Hamiltonian is sparsified. Initial conditions are also unaffected: ${\w{e}(t=0)}$ is set by interpolating ${\w{E}(t=0)}$ and scaling by $\mM_1$. To recover any subsequent field, it suffices to set ${\w{E}(t)\approx\mQ_1\cdot\w{e}(t)}$, such that $\mM_1^{-1}$ itself is altogether avoided.

For time integration, we decompose ${\tilde{H}_\text{d}=\tilde{H}_\w{e}+\tilde{H}_\w{a}}$ with ${\tilde{H}_\w{e}=\frac{\epsilon_0}{2}\w{e}^T\mQ_1\w{e}}$ and ${\tilde{H}_\w{a}=\frac{1}{2\mu_0}\w{a}^T\mC^T\tilde{\mM}_2\mC\w{a}}$. Each sub-Hamiltonian generates an exactly solvable shear map on $(\w{a},\w{e})$:
\begin{eqn}
\Phi^\w{e}_\tau:~&(\w{a},\w{e})\mapsto(\w{a}-\tau\,\mQ_1\,\w{e},~\w{e})\\
\Phi^\w{a}_\tau:~&(\w{a},\w{e})\mapsto(\w{a},~\w{e}+\tau\,c^2\mC^T\tilde{\mM}_2\mC\,\w{a}).
\label{shearMaps}
\end{eqn}
These compose via Strang splitting \cite{strang_construction_1968}, for example, to give the second-order map of time step $\tau$,
\begin{eqn}
\Phi_\tau=\Phi^\w{e}_{\tau/2}\circ\Phi^\w{a}_\tau\circ\Phi^\w{e}_{\tau/2},
\label{strangSplit}
\end{eqn}
though a GYM may employ any symplectic splitting composed of the same elemental maps of Eq.~(\ref{shearMaps}), including Lie--Trotter splitting \cite{trotter_product_1959} and higher-order compositions \cite{yoshida_construction_1990}. More generally, gauge-compatible splitting methods \cite{glasser_geometric_2020} provide a systematic framework for constructing such compositions while fully preserving gauge symmetry. We return to the role of gauge-compatible splittings in the formal definition of a GYM in Section~\ref{GYMdefnSubsect}.

\subsection{Mass-matrix independence of symplecticity\label{massMatrixIndep}}

Since the Poisson bracket $\{\cdot,\cdot\}_\text{d}$ is canonical and independent of the mass matrices, the flow generated by \emph{any} Hamiltonian with respect to this bracket is symplectic. The splitting integrator inherits this property, since each shear map is the exact flow of a sub-Hamiltonian. We state this formally and verify it by direct computation.

\begin{theorem}[Symplecticity of GYMs]\label{thm:symplecticity}
Let $\mQ_1\in\mR^{N_1\times N_1}$ and $\tilde{\mM}_2\in\mR^{N_2\times N_2}$ be arbitrary symmetric matrices. Then the one-step map $\Phi_\tau$ composed from the shear maps of Eq.~(\ref{shearMaps}) via any splitting scheme is symplectic with respect to the canonical symplectic form ${\mmOmega=\left[\begin{smallmatrix}\mZero&\mOne\\-\mOne&\mZero\end{smallmatrix}\right]}$ on $\mR^{2N_1}$.
\end{theorem}

\begin{proof}
Since the composition of symplectic maps is symplectic, it suffices to show that each shear map is individually symplectic. The Jacobian of $\Phi^\w{e}_\tau$ is
\begin{eqn}
\mJ_{\Phi^\w{e}_\tau}=\left[\begin{matrix}\mOne&-\tau\,\mQ_1\\\mZero&\mOne\end{matrix}\right].
\end{eqn}
Computing:
\begin{eqn}
(\mJ_{\Phi^\w{e}_\tau})^T\mmOmega\,\mJ_{\Phi^\w{e}_\tau}&=
\left[\begin{matrix}\mOne&\mZero\\-\tau\,\mQ_1^T&\mOne\end{matrix}\right]
\left[\begin{matrix}\mZero&\mOne\\-\mOne&\mZero\end{matrix}\right]
\left[\begin{matrix}\mOne&-\tau\,\mQ_1\\\mZero&\mOne\end{matrix}\right]=\left[\begin{matrix}\mZero&\mOne\\-\mOne&\tau(\mQ_1-\mQ_1^T)\end{matrix}\right]\\
&=\mmOmega,
\label{eqnUsesOmega}
\end{eqn}
where the last equality uses ${\mQ_1=\mQ_1^T}$, assumed WLOG. The calculation for $\Phi^\w{a}_\tau$ is analogous, and the result follows by composition.
\end{proof}

\begin{remark}
Theorem~\ref{thm:symplecticity} imposes no conditions on $\mQ_1$ or $\tilde{\mM}_2$ beyond symmetry ---which is assumed WLOG, since only their symmetric parts enter the EOM. Neither positive definiteness, sparsity, nor approximation quality affects symplecticity. Standard backward error analysis \cite{hairer_geometric_2006,benettin_hamiltonian_1994} then guarantees that $\Phi_\tau$ is the exact flow of a modified Hamiltonian $\tilde{H}_\epsilon = \tilde{H}_\text{d} + \mathcal{O}(\tau^p)$, where $p$ is the splitting order, up to exponentially small corrections.
\end{remark}

\subsection{Stability and the role of positive definiteness\label{stabilitySect}}

Although symplecticity places no constraints on the mass-matrix approximation, \emph{stability} does. The two properties are logically independent: a symplectic map can be unstable, and a stable map need not be symplectic. We shall find that to ensure the stable evolution of a splitting method composed of the symplectic maps of Eq.~(\ref{shearMaps}), two further conditions must be satisfied: (i) $\mQ_1$ and $\tilde{\mM}_2$ must be SPD; and (ii) the time interval $\tau$ must satisfy a Courant-Friedrichs-Lewy (CFL) condition.

We demonstrate this for the Strang splitting method of Eq.~(\ref{strangSplit}), though an analogous procedure may be followed for other splittings.

For notational convenience, we define the approximate curl-of-curl operator $\tilde{\mK}$ such that
\begin{eqn}
\tilde{\mK}=\mQ_1\mP_2~~\text{where}~~\mP_2=c^2\mC^T\tilde{\mM}_2\mC.
\end{eqn}
The evolution after $n$ Strang time steps is $(\mJ_{\Phi_\tau})^n$, where $\mJ_{\Phi_\tau}$ is the Jacobian of a single Strang time step $\tau$ as defined in Eq.~(\ref{strangSplit}):
\begin{eqn}
\mJ_{\Phi_\tau}=\left[\begin{matrix}\mOne-\frac{\tau^2}{2}\tilde{\mK} & \frac{\tau^3}{4}\tilde{\mK}\mQ_1-\tau\mQ_1\vspace{3pt}\\
\tau\mP_2 & \mOne-\frac{\tau^2}{2}\mP_2\mQ_1\end{matrix}\right].
\end{eqn}
Stability requires all eigenvalues of $\mJ_{\Phi_\tau}$ to lie on the unit circle. These eigenvalues are determined by the spectrum of $\tilde{\mK}$ \cite{leimkuhler_simulating_2004}. To see this, suppose we are given an eigenpair ${(\omega^2,\w{u}_{\omega})}$ of $\tilde{\mK}$. (Note ${\omega^2\geq0}$ because, for $\mQ_1$ and $\tilde{\mM}_2$ SPD, $\tilde{\mK}$ is similar to the symmetric positive semidefinite matrix ${\mQ_1^{1/2}\mP_2\mQ_1^{1/2}}$.) It is readily checked that
\begin{eqn}
V_{\w{u}_{\omega}}=\ws{span}\left\{\left(\begin{matrix}\w{u}_{\omega}\\\w{0}\end{matrix}\right),\left(\begin{matrix}\w{0}\\\mP_2\w{u}_{\omega}\end{matrix}\right)\right\}
\end{eqn}
is an invariant subspace of $\mJ_{\Phi_\tau}$. Indeed, the action of $\mJ_{\Phi_\tau}$ on ${V_{\w{u}_{\omega}}}$, expressed in the two dimensional basis above, is given by the ${2\times2}$ matrix
\begin{eqn}
\mJ_{\Phi_\tau}|_{V_{\w{u}_{\omega}}}=\left[\begin{matrix}1-\frac{\tau^2\omega^2}{2} & -\tau\omega^2\left(1-\frac{\tau^2\omega^2}{4}\right) \\
\tau & 1-\frac{\tau^2\omega^2}{2}\end{matrix}\right]
\end{eqn}
with eigenvalues
\begin{eqn}
\lambda^\pm = 1-\frac{z}{2}\pm\frac{1}{2}\sqrt{z(z-4)}
\label{eigenvaluePairs}
\end{eqn}
for ${z=\tau^2\omega^2}$, and satisfying ${\lambda^+\lambda^-=1}$. The reciprocal pair ${(\lambda^+,\lambda^-)}$ lie on the unit circle if and only if they are complex conjugates, which occurs when ${z(z-4)\leq0}$, or
\begin{eqn}
\tau\leq 2/\omega.
\label{singleEigCFLCond}
\end{eqn}
Eq.~(\ref{singleEigCFLCond}) is a CFL stability condition \cite{courant_uber_1928} for the Strang splitting evolution of Eq.~(\ref{discEOM}) for the particular eigenmode $\w{u}_\omega$ of $\tilde{\mK}$.

When $\mQ_1$ is symmetric but indefinite, $\tilde{\mK}$ can have negative eigenvalues $\omega^2<0$, producing exponential growth regardless of $\tau$---a structural instability that no timestep reduction can mitigate. Only when both $\mQ_1$ and $\tilde{\mM}_2$ are SPD is stability reducible to a CFL condition.

\begin{proposition}[Stability of GYMs]\label{prop:stability}
Let $\mQ_1$ and $\tilde{\mM}_2$ be symmetric positive definite. Then ${\tilde{\mK} = c^2\mQ_1\mC^T\tilde{\mM}_2\mC}$ has real nonnegative eigenvalues ${\{\omega_j^2\}}$, and the Strang splitting map is stable if and only if
\begin{eqn}
\tau \leq \frac{2}{\max_j\abs{\omega_j}}.
\label{fullCFLCond}
\end{eqn}
\end{proposition}

\begin{proof}
Eq.~(\ref{singleEigCFLCond}) must hold for each eigenvalue $\omega_j^2$ of $\tilde{\mK}$, so Eq.~(\ref{fullCFLCond}) follows.
\end{proof}

Theorem~\ref{thm:symplecticity} and Proposition~\ref{prop:stability} together establish that any SPD sparse approximation ${\mQ_1\approx\mM_1^{-1}}$ yields a GYM that is both symplectic and stable under the appropriate CFL condition. The SPD requirement is satisfied by all standard sparsification techniques. In practice, therefore, the choice of approximation is guided purely by accuracy and computational considerations, as we explore in Sections~\ref{YeeIsFEEC}--\ref{SparseMassMatSect}.

%

\subsection{Definition of a generalized Yee method\label{GYMdefnSubsect}}

We now formally define the class of algorithms studied in this paper.

\begin{definition}[Generalized Yee Method]\label{def:GYM}
A \emph{generalized Yee method} (GYM) is a finite element time-domain method for Maxwell's equations satisfying the following requirements: \textnormal{\textbf{(i)~de~Rham-conforming finite elements}}---the discrete fields are represented in a de~Rham-conforming finite element basis ${\{\gv{\Lambda}^p\}_{p=0}^n}$, with commuting projections as in Fig.~\ref{deRhamDiagram}; \textnormal{\textbf{(ii)~Sparse mass-matrix approximations}}---the exact inverse mass matrix $\mM_1^{-1}$ is replaced by a sparse SPD approximation ${\mQ_1\approx\mM_1^{-1}}$, and $\mM_2$ is optionally replaced by a sparse SPD ${\tilde{\mM}_2\approx\mM_2}$, yielding the approximate Hamiltonian $\tilde{H}_\text{d}$ of Eq.~(\ref{discTildeHamilSystem}); \textnormal{\textbf{(iii)~Gauge-compatible splitting}}---time evolution proceeds by composition of the shear maps of Eq.~(\ref{shearMaps}) via a gauge-compatible splitting method \cite{glasser_gauge-compatible_2022,glasser_geometric_2020}.
\end{definition}

%
%
%

\begin{table*}[t!]
\renewcommand{\arraystretch}{1.1}
\centering
\begin{tabular*}{.97875\textwidth}{|c|c|c|}
\hline
 & \textbf{Yee's Method} & \textbf{Generalized Yee Methods}\\
 \hline
Mesh & cubical & simplicial or cubical\\
\hline
Finite elements & Whitney forms & any${}^*$ de~Rham conforming\\
\hline
$\mQ_1$ and $\tilde{\mM}_2$  & diagonal (lumped) & any sparse SPD approximations\\
\hline
Splitting method & Strang & any gauge-compatible splitting\\
\hline
\end{tabular*}
\caption{Yee's method and its generalization to GYMs. Each row corresponds to a design choice generalized from Yee's method to GYMs, as formalized in Definition~\ref{def:GYM}. ${}^*$For symplectic particle-coupling (as in PIC simulations), the finite elements must be at least $C^1$ \cite{barham_diagnosing_2025}; see Section~\ref{PICsect}.}
\label{YeeVsGYM}
\end{table*}

Each of these requirements generalizes a specific aspect of Yee's method, as summarized in Table~\ref{YeeVsGYM}. Requirement~(i) ensures the discrete fields inherit the topological properties of the de~Rham complex (Fig.~\ref{deRhamDiagram}): divergence-free magnetic fields (${\mmD\mC=\mZero}$) and gauge symmetry, independent of any mass-matrix approximation. Requirement~(ii), by Theorem~\ref{thm:symplecticity}, preserves symplecticity for any symmetric $\mQ_1$ and $\tilde{\mM}_2$; the SPD condition ensures stability (Proposition~\ref{prop:stability}). The sparsity of $\mQ_1$ and $\tilde{\mM}_2$ restricts each shear map to local data communication, making GYMs suitable for massively parallel computation.

Requirement~(iii) specifies that the splitting be gauge-compatible in the sense of \cite{glasser_gauge-compatible_2022, glasser_geometric_2020}. A splitting method for a GYM is gauge-compatible when (1) each sub-Hamiltonian is individually invariant under the discrete gauge transformation ${\w{a}\mapsto\w{a}+\mG\w{s}}$, and (2) each subsystem is solved exactly. When both conditions hold, the momentum map associated with the gauge symmetry---in the case of Maxwell's equations, the momentum map is simply Gauss' law---is exactly preserved at every timestep, guaranteeing exact charge conservation after time discretization. For the source-free system, gauge compatibility is automatic: $\tilde{H}_\w{a}$ depends on $\w{a}$ only through ${\mC\w{a}}$, and ${\mC\mG=\mZero}$ renders it gauge-invariant; $\tilde{H}_\w{e}$ is independent of $\w{a}$ entirely. Both sub-Hamiltonians generate linear flows that are solved exactly by the shear maps of Eq.~(\ref{shearMaps}). We see again that the mass-matrix approximations $\mQ_1$ and $\tilde{\mM}_2$ modify only the metric structure entering the Hamiltonian, while the differential operators $\mG$ and $\mC$---which underpin gauge invariance via ${\mC\mG=\mZero}$---remain exact. Gauge compatibility is therefore preserved by any GYM. This symmetry constraint will be revisited in our discussion of particle-coupling in Section~\ref{PICsect}.

In summary, GYMs form a family of algorithms parameterized by four choices: the mesh, the finite element basis, the SPD mass-matrix approximation, and the splitting scheme. The symplecticity and topological preservation of the method are guaranteed by the structure of Definition~\ref{def:GYM}, independent of the specific choices made within each category. The practitioner is therefore free to optimize each choice based on accuracy and computational considerations alone.

\section{Classifying Yee's Method as a GYM\label{YeeIsFEEC}}

We now show that Yee's method is recovered as the simplest GYM: generalized Whitney forms on a cubical mesh, diagonal mass matrices, and Strang splitting. After establishing these ingredients, we verify that the resulting time advance reproduces Yee's update equations.

\subsection{Whitney forms on a cubical mesh}

We consider a cubical lattice ${\mT_h\subset\mR^3}$ with lattice spacings ${\{\Delta_x,\Delta_y,\Delta_z\}}$ and choose the \emph{generalized Whitney forms} ${Q_1^-\Lambda^p(\mT_h)}$ \cite{arnold_periodic_nodate} as our FEEC basis, a family of piecewise polynomial finite elements also defined in \cite{lohi_whitney_2021} Example~5.2.

A generalized Whitney $p\text{-form}$ is defined by its 1-to-1 correspondence with a $p$-dimensional feature of the mesh. The generalized Whitney $p\text{-form}$ ${\mW_{\sigma^p}\in Q_1^-\Lambda^p(\mT_h)}$ associated to a given $p$-face ${\sigma^p\subset\mT_h\subset\mR^n}$ is defined by
\begin{eqn}
\frac{1}{\abs{\nu^p}}\int_{\nu^p}\mW_{\sigma^p} =
\begin{cases}
1 & \nu^p=\sigma^p\\
0&\nu^p\neq\sigma^p,
\end{cases}
\label{WhitneyFormDefn}
\end{eqn}
so that generalized Whitney $p\text{-form}$s are \emph{dual} to $p$-faces of $\mT_h$ via integration. Here, ${\abs{\nu^p}}$ denotes the $p$-volume of $\nu^p$ (1, length, area, or volume for ${p=0,1,2,3}$, respectively). This normalization defines fields as they are typically represented in Yee's method.

\begin{figure}[b!]
\centering
\includegraphics[width=0.48\columnwidth]{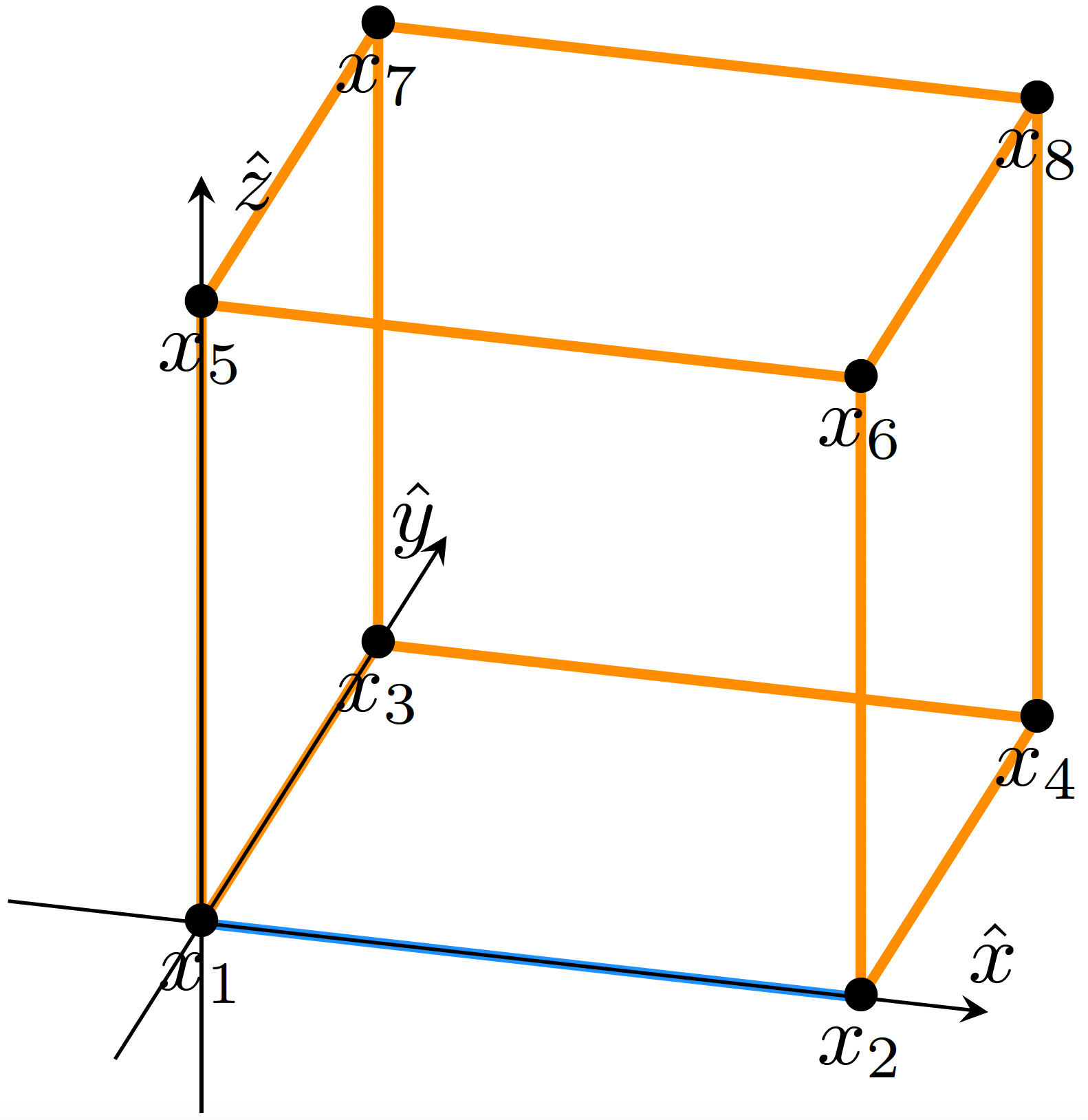}
\hfill
\includegraphics[width=0.48\columnwidth]{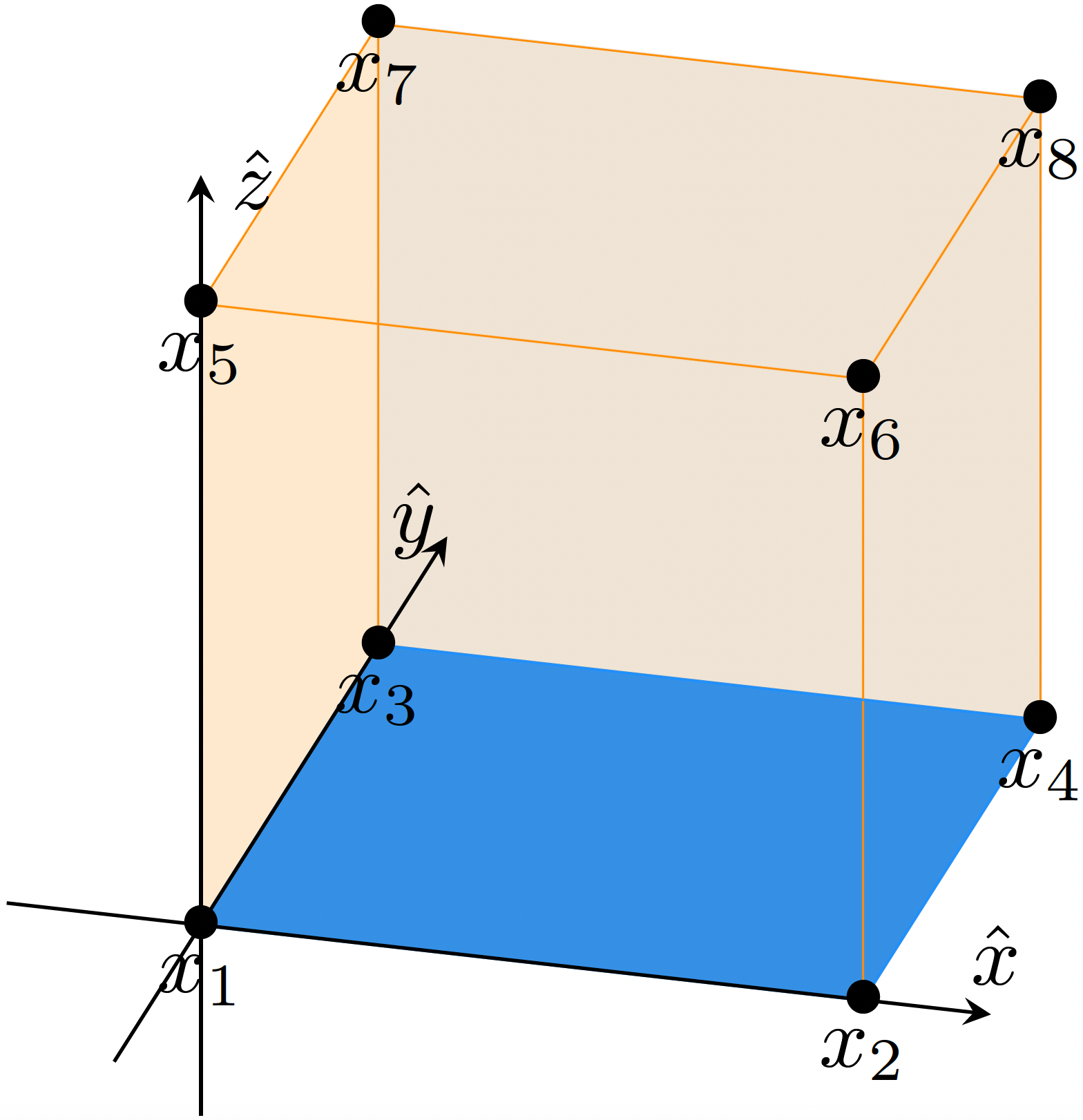}
\caption{The generalized Whitney $\text{1-form}$ $\mW_{x_1x_2}$ (left) and $\text{2-form}$ $\mW_{x_1x_2x_4x_3}$ (right). $\mW_{x_1x_2}$ evaluates to $\abs{x_1x_2}$ when integrated along the edge ${x_1x_2}$ (blue) and vanishes on all other edges (orange). $\mW_{x_1x_2x_4x_3}$ likewise yields $\abs{x_1x_2x_4x_3}$ on the blue face and vanishes on all others.}
\label{GenWhitneyCubes}
\end{figure}

A concrete example is depicted in Fig.~\ref{GenWhitneyCubes}. The $\text{1-form}$ ${\mW_{x_1x_2}\in Q_1^-\Lambda^1(\mT_h)}$ and the $\text{2-form}$ ${\mW_{x_1x_2x_4x_3}\in Q_1^-\Lambda^2(\mT_h)}$ are given by
\begin{eqn}
\mW_{x_1x_2}&=\left(1-\tfrac{y}{\Delta_y}\right)\left(1-\tfrac{z}{\Delta_z}\right)\md x\\
\mW_{x_1x_2x_4x_3}&=\left(1-\tfrac{z}{\Delta_z}\right)\md x\wedge\md y,
\label{whitneyEx12}
\end{eqn}
and the Whitney $\text{2-form}$ associated to the face ${x_1x_5x_6x_2}$ is
\begin{eqn}
\mW_{x_1x_5x_6x_2}&=\left(1-\tfrac{y}{\Delta_y}\right)\md z\wedge\md x.
\label{whitneyEx12_a}
\end{eqn}
As required by Eq.~(\ref{WhitneyFormDefn}), $\mW_{x_1x_2}$ evaluates to $\abs{x_1x_2}$ when integrated along the edge $x_1x_2$ and vanishes on all other edges; $\mW_{x_1x_2x_4x_3}$ evaluates to $\abs{x_1x_2x_4x_3}$ on its associated face and vanishes on all others.

\subsection{The curl operator as a finite difference}

The curl matrix ${\mC:\mR^{N_1}\rightarrow\mR^{N_2}}$ computes $\md$ for $\text{1-form}$s, mapping coefficients in the ${Q_1^-\Lambda^1(\mT_h)}$ basis to those in ${Q_1^-\Lambda^2(\mT_h)}$. Because generalized Whitney forms are in 1-to-1 correspondence with mesh features, $\mC$ admits a simple interpretation: directly computing $\md\mW_{x_1x_2}$ from Eqs.~(\ref{whitneyEx12}--\ref{whitneyEx12_a}) gives
\begin{eqn}
\md\mW_{x_1x_2}&=\frac{1}{\Delta_y}\mW_{x_1x_2x_4x_3}-\frac{1}{\Delta_z}\mW_{x_1x_5x_6x_2}.
\label{dCurlExample}
\end{eqn}
A Whitney $\text{2-form}$ $\mW_{\sigma^2}$ appears on the right-hand side if and only if the face $\sigma^2$ contains the edge $x_1x_2$ on its boundary. Its coefficient ${\pm1/\Delta_{x^\mu}}$ is determined by the dimension in which $\sigma^2$ extends $x_1x_2$, with sign set by the relative orientation.

\begin{figure}[t!]
\includegraphics[width=\linewidth]{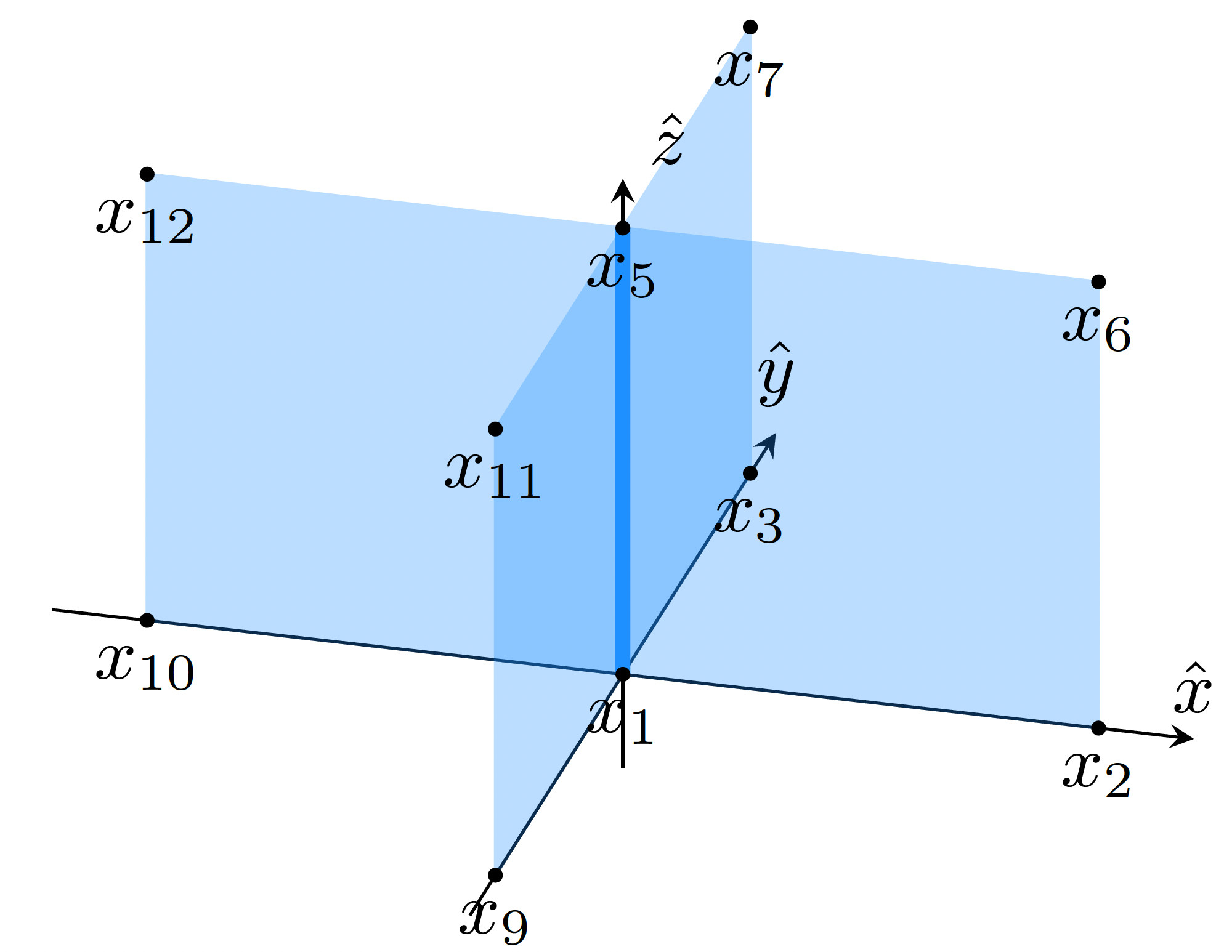}
\caption{Degrees of freedom involved in the transposed curl operator $\mC^T$ for generalized Whitney forms on a cubic mesh.}
\label{CTranspGraphic}
\end{figure}

Thus $\mC$ acts on generalized Whitney form coefficients as a finite difference operator. Denoting by ${\w{a}_{\sigma^1}}$and ${\w{b}_{\sigma^2}}$ the respective coefficients of $\mW_{\sigma^1}$ and $\mW_{\sigma^2}$ in ${\w{b}=\mC\w{a}}$, we find (with reference to Fig.~\ref{GenWhitneyCubes})
\begin{eqn}
\w{b}_{x_1x_2x_4x_3}&=(\mC\w{a})_{x_1x_2x_4x_3}=\frac{1}{\Delta_x}\Big(\w{a}_{x_2x_4}-\w{a}_{x_1x_3}\Big)-\frac{1}{\Delta_y}\Big(\w{a}_{x_3x_4}-\w{a}_{x_1x_2}\Big).
\label{dIsFiniteDiff}
\end{eqn}
The transpose $\mC^T$ yields an analogous finite difference between faces adjoining a given edge. Using Fig.~\ref{CTranspGraphic} as a reference,
\begin{eqn}
(\mC^T\w{b})_{x_1x_5}=\frac{1}{\Delta_x}\Big(\w{b}&_{x_1x_5x_6x_2}-\w{b}_{x_{10}x_{12}x_5x_1}\Big)-\frac{1}{\Delta_y}\Big(\w{b}_{x_1x_3x_7x_5}-\w{b}_{x_{9}x_{1}x_5x_{11}}\Big).
\label{dTIsFiniteDiff}
\end{eqn}

\subsection{Diagonal mass matrices}

For Whitney forms on a uniform cubical mesh, the mass matrix entries are (away from the boundary)
\begin{eqn}
(\mM_1)_{\sigma^1,\nu^1}&=\Delta_V\cdot
\begin{cases}
4/9&\sigma^1=\nu^1\\
1/9&\sigma^1\parallel\nu^1\text{ and }\exists~\nu^2\supset\{\sigma^1,\nu^1\}\\
1/36&\sigma^1\parallel\nu^1\text{ and }\exists~\nu^3\supset\{\sigma^1,\nu^1\}
\end{cases}\\
(\mM_2)_{\sigma^2,\nu^2}&=\Delta_V\cdot
\begin{cases}
2/3&\sigma^2=\nu^2\\
1/6&\sigma^2\parallel\nu^2\text{ and }\exists~\nu^3\supset\{\sigma^2,\nu^2\}
\end{cases}
\label{whitMassMatrices}
\end{eqn}
where ${\Delta_V=\Delta_x\Delta_y\Delta_z}$ denotes a cell volume. Each row of $\mM_1$ and $\mM_2$ sums to $\Delta_V$, motivating the lumped approximation
\begin{eqn}
\mM_1^Y&=\Delta_V\cdot\mOne_{N_1\times N_1}\,,\qquad
\mM_2^Y=\Delta_V\cdot\mOne_{N_2\times N_2}.
\label{YeeMassMatrices}
\end{eqn}
By Theorem~\ref{thm:symplecticity}, this lumping modifies the Hamiltonian but leaves the symplectic structure undisturbed.

\subsection{Recovering Yee's update equations}

We now verify that the Strang splitting of Eq.~(\ref{strangSplit}), applied to Eq.~(\ref{preYeeFEECEqns}) with the lumped mass matrices of Eq.~(\ref{YeeMassMatrices}), reproduces Yee's method. In terms of ${\w{b}=\mC\w{a}}$, the shear maps become
\begin{eqn}
\Phi^\w{e}_\tau:~&(\w{b},\w{e})\mapsto\big(\w{b}-(\tau/\Delta_V)\,\mC\,\w{e},~\w{e}\big)\\
\Phi^\w{b}_\tau:~&(\w{b},\w{e})\mapsto\big(\w{b},~\w{e}+\tau\,\Delta_V c^2\,\mC^T\,\w{b}\big).
\label{shearMapsBE}
\end{eqn}

An initial half-step $\Phi^\w{e}_{\tau/2}$ offsets $\w{b}$ by a half timestep, establishing the staggered initial data characteristic of Yee's method. Thereafter, each full step consists of $\Phi^\w{b}_\tau$ followed by $\Phi^\w{e}_\tau$, yielding the leapfrog updates:
\begin{eqn}
\Phi^\w{b}_\tau:~~\w{e}\left[t_{n+1}\right]&=\w{e}\left[t_n\right]+\tau\,\Delta_V c^2\,\mC^T\w{b}\!\left[t_{n+1/2}\right]\\
\Phi^\w{e}_\tau:~~\w{b}\left[t_{n+3/2}\right]&=\w{b}\!\left[t_{n+1/2}\right]-(\tau/\Delta_V)\,\mC\,\w{e}\!\left[t_{n+1}\right].
\label{YeeLeapfrog}
\end{eqn}
Since $\mC$ and $\mC^T$ act as finite differences on the cubical Whitney form basis (Eqs.~(\ref{dIsFiniteDiff})--(\ref{dTIsFiniteDiff})), these are precisely the Yee update equations. For example, Eq.~(\ref{YeeLeapfrog}a) gives
\begin{eqn}
&\frac{1}{c^2 \tau}\left(E_x^{n+1}\left[i+\tfrac{1}{2},j,k\right]-E_x^n\left[i+\tfrac{1}{2},j,k\right]\right)\\
&=\frac{1}{\Delta_y}\left(B_z^{n+\frac{1}{2}}\left[i+\tfrac{1}{2},j+\tfrac{1}{2},k\right]-B_z^{n+\frac{1}{2}}\left[i+\tfrac{1}{2},j-\tfrac{1}{2},k\right]\right)\\
&-\frac{1}{\Delta_z}\left(B_y^{n+\frac{1}{2}}\left[i+\tfrac{1}{2},j,k+\tfrac{1}{2}\right]-B_y^{n+\frac{1}{2}}\left[i+\tfrac{1}{2},j,k-\tfrac{1}{2}\right]\right),
\label{Yee1}
\end{eqn}
and Eq.~(\ref{YeeLeapfrog}b) gives
\begin{eqn}
\frac{1}{\tau}&\left(B_x^{n+\frac{1}{2}}\left[i,j+\tfrac{1}{2},k+\tfrac{1}{2}\right]-B_x^{n-\frac{1}{2}}\left[i,j+\tfrac{1}{2},k+\tfrac{1}{2}\right]\right)\\
&=\frac{1}{\Delta_z}\left(E_y^n\left[i,j+\tfrac{1}{2},k+1\right]-E_y^n\left[i,j+\tfrac{1}{2},k\right]\right)\\
&-\frac{1}{\Delta_y}\left(E_z^n\left[i,j+1,k+\tfrac{1}{2}\right]-E_z^n\left[i,j,k+\tfrac{1}{2}\right]\right),
\label{Yee2}
\end{eqn}
up to the constant factor $\Delta_V$, which reflects the normalization convention of Eq.~(\ref{WhitneyFormDefn}) and affects only the units of $\w{b}$, not the dynamics. Yee's algorithm, therefore, is the GYM defined by Whitney forms on a cubical mesh, lumped mass matrices, and Strang splitting, as summarized in Table~\ref{YeeVsGYM}.

We have further confirmed this equivalence numerically: implementing the full GYM pipeline with ${\mQ_1=(\Delta_V)^{-1}\mOne}$ and ${\tilde{\mM}_2=\Delta_V\mOne}$ on a cubical mesh reproduces the standard Yee leapfrog updates to machine precision, with ${\norm{\w{e}^\text{GYM}(t)-\w{e}^\text{Yee}(t)}/\norm{\w{e}^\text{Yee}(t)}<10^{-15}}$ after a full time step. More such numerical results will be presented in Section~\ref{NumResult}.

\section{Sparse mass matrix approximations for GYMs\label{SparseMassMatSect}}

While FEEC mass matrices $\mM_1$ and $\mM_2$ are themselves sparse, the EOM of Eq.~(\ref{preYeeFEECEqns}) require the typically dense inverse $\mM_1^{-1}$. This density couples DOFs globally, requiring communication at every timestep between all computational nodes and spoiling scalability. By Theorem~\ref{thm:symplecticity}, any sparse SPD approximation $\mQ_1\approx\mM_1^{-1}$ yields a symplectic and stable GYM, so one is free to choose a sparsification strategy optimizing accuracy and computational cost. Given this flexibility, we briefly survey existing approaches in the literature, before introducing the SPAI-OP method.

\paragraph{Mass lumping} As in Yee's method (see Section~\ref{YeeIsFEEC}), the simplest approach replaces $\mM_p$ with a lumped diagonal matrix $\mM_p^L$ whose $i^\text{th}$ entry is the sum of the $i^\text{th}$ row of $\mM_p$, i.e.\ ${(\mM_p^L)_{ii} = \sum_j (\mM_p)_{ij}}$ \cite{cohen_gauss_1998}. For the lowest-order Whitney forms, this incurs errors consistent with the convergence rate of a first-order finite element space, so no accuracy is sacrificed. For higher-order elements, however, naive row-sum lumping discards off-diagonal elements, effectively reducing the mass matrix to a low-order quadrature approximation regardless of polynomial degree.

\paragraph{Spectral elements} Spectral element methods overcome this limitation through a mutually compatible design of basis functions and Gauss--Lobatto--Legendre (GLL) quadrature \cite{duczek_mass_2019}. However, GLL spectral elements are $C^0$ across cell boundaries, which is insufficient for PIC methods requiring $C^1$ field smoothness \cite{barham_diagnosing_2025} (see Section~\ref{PICsect}).

\paragraph{Broken finite element spaces}
Broken finite element spaces \cite{campos_pinto_gauss-compatible_2016,kapidani_arbitrary-order_2021,guclu_broken_2023} restrict each basis function to a single cell, rendering $\mM_p$ block-diagonal and exactly (and sparsely) invertible. However, continuity must be restored through lifting operators or penalty terms, multiplying the number of degrees of freedom and requiring care to maintain the discrete de~Rham complex of Fig.~\ref{deRhamDiagram} \cite{campos_pinto_gauss-compatible_2016}.

\paragraph{Sparse approximate inverses}
A more direct strategy computes a sparse approximation ${\mQ_p \approx \mM_p^{-1}}$ over a prescribed sparsity pattern. The pattern may be chosen by thresholding the entries of $\mM_p^{-1}$ below a tolerance \cite{bo_he_sparse_2006,he_differential_2007}, or by retaining couplings up to a topological distance, as identified by the sparsity pattern of a power $(\mM_p)^k$ of $\mM_p$. Given the sparsity pattern, the entries of $\mQ_p$ are determined by minimizing $\norm{\mM_p\mQ_p - \mOne}_F$. The Frobenius norm separates column-wise \cite{huckle_approximate_1998}:
\begin{eqn}
\norm{\mM_p\mQ_p-\mOne}_F^2=\sum\limits_{\ell=1}^{N_p}\norm{\mM_p\w{q}_\ell-\mathbf{1}_\ell}^2,
\label{FrobeniusDecomp}
\end{eqn}
where $\w{q}_\ell$ is the $\ell^{\text{th}}$ column of ${\mQ_p}$ and $\mathbf{1}_\ell$ the $\ell^{\text{th}}$ standard basis vector. Let ${\mI_\ell\subset\{1,\dots,N_1\}}$ index the nonzero entries permitted in $\w{q}_\ell$ by the prescribed sparsity pattern. The column-wise minimization of Eq.~(\ref{FrobeniusDecomp}) admits the closed-form solution
\begin{eqn}
\w{q}_\ell(\mI_\ell)=\left[\mM_p(:,\mI_\ell)^T\mM_p(:,\mI_\ell)\right]^{-1}\mM_p(:,\mI_\ell)^T\mathbf{1}_\ell,
\label{NormalEqSolve}
\end{eqn}
a small ${\abs{\mI_\ell}\times\abs{\mI_\ell}}$ system solved directly---an embarrassingly parallel computation \cite{kim_parallel_2011,teixeira_differential_2013,grote_spai_1997,huckle_approximate_1998}. Since each column is optimized independently, the resulting $\mQ_p$ is not inherently symmetric. Symmetry---required by Proposition~\ref{prop:stability} for stability---may be enforced by setting ${\mQ_p\leftarrow\frac{1}{2}(\mQ_p+\mQ_p^T)}$ at the cost of a small (typically ${\lesssim10\%}$) increase in the Frobenius error of Eq.~(\ref{FrobeniusDecomp}). This added symmetrization error can be partially offset by gradient descent through symmetric matrix space using techniques described in Section~\ref{SPAIOPsubsect}.

\subsection{SPAI-OP: Operator-probed sparse approximate inverses\label{SPAIOPsubsect}}

Standard SPAI optimizes $\mQ_1$ with respect to inversion alone, and so distributes approximation error uniformly across all modes of the identity operator. However, in the GYM context, $\mQ_1$ appears within the approximate curl-of-curl operator ${\tilde{\mK}=\mQ_1\mP_2}$, where ${\mP_2\equiv c^2\mC^T\tilde{\mM}_2\mC}$. When ${\mQ_1=\mM_1^{-1}}$ is exact, $\tilde{\mK}$ reduces to the exact discrete curl-of-curl operator ${\mK\equiv\mM_1^{-1}\mP_2}$, which inherits the full accuracy of the finite element discretization. For a sparse $\mQ_1$, the sparsification error on a mode $\w{u}$ is
\begin{eqn}
(\tilde{\mK}-\mK)\w{u} = (\mQ_1-\mM_1^{-1})\mP_2\w{u} \equiv \delta\mQ_1\,\mP_2\w{u},
\label{sparsificationError}
\end{eqn}
where ${\delta\mQ_1\equiv\mQ_1-\mM_1^{-1}}$. This error is small precisely when $\mQ_1$ accurately approximates the action of $\mM_1^{-1}$ on $\mP_2\w{u}$---a geometric condition reflecting the wave-propagation structure of the GYM rather than the abstract matrix inversion problem.

In some applications, a few physically distinguished wave modes $\w{v}_k$ dominate the accuracy requirements---for example, eigenmodes of $\mK$ within a frequency band of interest, particular initial conditions, or wave packets relevant to a beam-driven instability or resonance. Concentrating accuracy on such modes (even at the cost of uniform spectral accuracy) is desirable.  Inspired by the \emph{modified SPAI} (MSPAI) framework of Huckle and Kallischko \cite{huckle_frobenius_2007}, who introduced probing constraints to improve sparse approximate inverses on targeted subspaces for iterative solver preconditioning, we develop an \emph{operator-probed SPAI} (SPAI-OP) method for the wave-propagation context of GYMs.  SPAI-OP augments the standard SPAI objective with soft constraints that penalize the sparsification error of Eq.~(\ref{sparsificationError}) on the target modes:
\begin{eqn}
\mO_\text{SPAI-OP}(\mQ_1) = \norm{\mM_1\mQ_1-\mOne}_F^2
  + \lambda\sum_{k=1}^{n_\text{probe}}\norm{(\mM_1\mQ_1-\mOne)\mP_2\w{v}_k}^2,
\label{SPAIOPobjective}
\end{eqn}
where $\lambda>0$ balances the baseline Frobenius fit with the probing penalty.  The probing vectors $\mP_2\w{v}_k$ represent the action of the discrete curl-of-curl operator on the target modes $\w{v}_k$: their presence in the objective derives precisely from the geometric context of the wave operator, distinguishing SPAI-OP from the MSPAI methods that inspired it.

\paragraph{Symmetry-constrained Sylvester formulation} We seek a symmetric $\mQ_1$ supported on a prescribed symmetric sparsity pattern ${S\subset[N_1]\times[N_1]}$, where ${[N_1]\equiv\{1,\dots,N_1\}}$, satisfying the SPD requirement of Proposition~\ref{prop:stability}. Because $(\mQ_1)_{ij}$ and $(\mQ_1)_{ji}$ share a single degree of freedom, the stationarity condition for minimizing $\mO_\text{SPAI-OP}$ over symmetric $\mQ_1$ is
\begin{eqn}
\frac{\partial\mO_\text{SPAI-OP}}{\partial (\mQ_1)_{ij}}
  + \frac{\partial\mO_\text{SPAI-OP}}{\partial (\mQ_1)_{ji}} = 0,
  \qquad \forall\,(i,j)\in S.
\label{symmetricOptimality}
\end{eqn}
Defining ${\mU\equiv\mM_1^2}$, the probing matrix ${\mmW \equiv \mP_2\mV}$ with ${\mV=[\w{v}_1,\dots,\w{v}_{n_\text{probe}}]}$, and
\begin{eqn}
\hat{\mmT} \equiv \mOne + \lambda\,\mmW\mmW^T, \qquad
\hat{\mR} \equiv \mM_1\hat{\mmT} + \hat{\mmT}\mM_1,
\label{ThatDefinition}
\end{eqn}
Eqs.~(\ref{SPAIOPobjective}) and (\ref{symmetricOptimality}) yield the generalized Sylvester equation
\begin{eqn}
(\mU\mQ_1\hat{\mmT} + \hat{\mmT}\mQ_1\mU)_{ij} = \hat{\mR}_{ij},
  \qquad \forall\,(i,j)\in S.
\label{SPAIOPSylvester}
\end{eqn}
Setting $\lambda=0$ recovers symmetry-constrained SPAI ($\hat{\mmT}=\mOne$, $\hat{\mR}=2\mM_1$), whose solution is the \emph{optimal} symmetric $\mQ_1$ with respect to the Frobenius objective $\norm{\mM_1\mQ_1-\mOne}_F^2$ over the sparsity pattern $S$---a smaller error than the post-symmetrized column-wise SPAI of Eq.~(\ref{NormalEqSolve}), which merely projects the unconstrained minimizer onto the symmetric subspace. Increasing $\lambda$ concentrates accuracy on the probed modes.

\paragraph{Vectorized form and automatic symmetry of the solution} Let $S$ contain $n_q$ entries indexed by row-column pairs $(i_r(k),i_c(k))$ for ${k=1,\dots,n_q}$.  Setting ${\w{q}_k=(\mQ_1)_{i_r(k),i_c(k)}}$ and ${\w{r}_j=\hat{\mR}_{i_r(j),i_c(j)}}$, Eq.~(\ref{SPAIOPSylvester}) takes the equivalent vectorized form ${\mA\w{q}=\w{r}}$ with
\begin{eqn}
\mA_{jk} = \mU_{i_r(j),i_r(k)}\,\hat{\mmT}_{i_c(k),i_c(j)}
  + \hat{\mmT}_{i_r(j),i_r(k)}\,\mU_{i_c(k),i_c(j)}.
\label{SPAIOPcoefficient}
\end{eqn}
$\mA$ is a principal submatrix of the Kronecker sum ${\mU\otimes\hat{\mmT}+\hat{\mmT}\otimes\mU}$, hence symmetric positive semidefinite (since $\mU$ and $\hat{\mmT}$ are). A small Tikhonov regularization ${\mA\rightarrow\mA+\epsilon\mOne}$ ensures strict positive definiteness.

The solution of $\mA\w{q}=\w{r}$ yields a symmetric $\mQ_1$ whenever $S$ is a symmetric sparsity pattern. To see this, let ${\mPi\in\mR^{n_q\times n_q}}$ be the permutation matrix that swaps each sparsity index $k$ (corresponding to entry $(i_r(k),i_c(k))$) with the index associated with the transposed entry $(i_c(k),i_r(k))$.  The two Kronecker terms in Eq.~(\ref{SPAIOPcoefficient}) exchange under this permutation, so ${\mPi\mA=\mA\mPi}$; likewise, since $\hat{\mR}$ is symmetric, ${\mPi\w{r}=\w{r}}$.  From $\mA\w{q}=\w{r}$ it follows that $\mA(\mPi\w{q})=\mPi\mA\w{q}=\mPi\w{r}=\w{r}$, so $\mPi\w{q}$ is also a solution.  Since $\mA$ is strictly positive definite, the solution is unique, and ${\mPi\w{q}=\w{q}}$---equivalently, ${(\mQ_1)_{ij}=(\mQ_1)_{ji}}$.

The symmetry-constrained Sylvester formulation therefore confers two benefits over the column-wise SPAI of Eq.~(\ref{NormalEqSolve}): it produces a symmetric solution directly (without the post-symmetrization error discussed in the preceding subsection), and it accommodates the geometric probing constraints of Eq.~(\ref{SPAIOPobjective}) through a simple modification of $\hat{\mmT}$.

\paragraph{Matrix-free PCG solver} The coefficient matrix ${\mA\in\mR^{n_q\times n_q}}$ is too large to form explicitly for practical problems. The product $\mA\w{q}$ is nevertheless readily computed matrix-free: assemble the sparse matrix $\mQ_1$ from $\w{q}$ (on the sparsity pattern $S$), compute ${\mS=\mU\mQ_1\hat{\mmT}+\hat{\mmT}\mQ_1\mU}$, and extract ${(\mA\w{q})_j=\mS_{i_r(j),i_c(j)}}$. The symmetric positive definiteness of $\mA$ then allows preconditioned conjugate gradient (PCG) iteration. A diagonal preconditioner using the diagonal entries of $\mA$,
\begin{eqn}
\ws{prec}_j = \mU_{i_r(j),i_r(j)}\hat{\mmT}_{i_c(j),i_c(j)}
  + \hat{\mmT}_{i_r(j),i_r(j)}\mU_{i_c(j),i_c(j)} + \epsilon,
\label{diagonalPreconditioner}
\end{eqn}
is effective in practice: for ${\lambda=0}$ (pure symmetric SPAI), PCG converges to machine precision in $2$ iterations on structured periodic grids, essentially independently of mesh size.  For ${\lambda>0}$ the conditioning of the Sylvester operator deteriorates as the rank-$n_\text{probe}$ term $\lambda\mmW\mmW^T$ grows, and the iteration count grows commensurately; at the tolerance-optimal values $\lambda^*$ adopted in Section~\ref{NumResult} this remains modest (${\sim}10^2$ iterations on our structured test grids), and the Sylvester solve is in any case a one-time setup cost, its result amortized over the full time-domain simulation.  The efficient column-wise SPAI of Eq.~(\ref{NormalEqSolve}), followed by post-symmetrization, provides a natural warm-start for the PCG iteration of Eq.~(\ref{SPAIOPSylvester}): since the column-wise SPAI optimizes the same unweighted Frobenius objective, its symmetrized result closely approximates the ${\lambda=0}$ solution.

\paragraph{Choice of probing vectors and practical considerations} The probing weight $\lambda$ is a single user-facing parameter. Setting ${\lambda=0}$ recovers standard symmetric SPAI; increasing $\lambda$ concentrates accuracy on the probed modes at the expense of the overall Frobenius fit. The sensitivity analysis in Section~\ref{lambdaScan} shows that intermediate values of $\lambda$ substantially improve probed-mode accuracy while maintaining overall error near the SPAI baseline.

The sparsity pattern $S$ is typically chosen as ${S(\mQ_1) = S(\mM_1^k)}$ for some ${k=0,1,2,\dots}$, corresponding to stencils of increasing radius. Larger stencils yield more accurate approximations at higher per-timestep cost. Positive definiteness of $\mQ_1$, required for stability (Proposition~\ref{prop:stability}), is typically inherited from the SPD structure of $\mM_1$ when the sparsity pattern includes the diagonal and is not too coarse; in practice, positive definiteness should be verified after construction and enforced if necessary by a small diagonal shift.

The choice of probing vectors $\w{v}_k$ is application-dependent. For dispersion optimization, the eigenmodes of ${\mK=\mM_1^{-1}\mP_2}$ at physically relevant frequencies are natural targets, since these correspond to the resolved wave modes whose accurate propagation is most critical. For particle-in-cell simulations, modes resonant with beam-driven instabilities may be prioritized (see Section~\ref{PICsect}). Because the target modes need only be determined once at setup time---and only approximately, since the probing constraint is soft---the cost of computing them is typically a small fraction of total GYM setup time.

\section{GYMs for PIC\label{PICsect}}

We now extend the GYM framework to structure-preserving electromagnetic particle-in-cell (PIC) methods. The concern is subtle: while symplecticity of the field-only Hamiltonian flow is guaranteed by Theorem~\ref{thm:symplecticity} irrespective of the mass-matrix approximation, the coupling of particles to the fields through spatial \emph{interpolation} introduces a distinct requirement---namely, that the GYM Hamiltonian vector field be spatially continuous. We now review the argument of Barham and Burby \cite{barham_diagnosing_2025}, which identifies this continuity as a necessary condition for symplecticity over particle trajectories, and translates it into a pointwise $C^1$ smoothness requirement on the 1-form finite element basis.

\subsection{The Poincar\'e integral invariant and discrete smoothness}

A Hamiltonian flow $\Phi_t$ preserves the canonical symplectic form, ${\Phi_t^*\omega_c = \omega_c}$, where ${\omega_c = d\w{q}\wedge d\w{p}}$. Writing ${\omega_c = -d\vartheta}$ with ${\vartheta = \w{p}\cdot d\w{q}}$ the Liouville 1-form, one obtains the first Poincar\'e integral invariant: for any closed phase-space loop $\gamma = \partial D$,
\begin{eqn}
I(t) \;=\; \oint_{\Phi_t \circ \gamma} \vartheta \;=\; -\int_{\Phi_t(D)}\omega_c \;=\; -\int_D \Phi_t^*\omega_c \;=\; -\int_D \omega_c \;=\; \oint_\gamma \vartheta \;=\; I(0),
\label{PoincareInvariant}
\end{eqn}
where the second and penultimate equalities follow from Stokes' theorem, the third from the change-of-variables formula, and the fourth from the preservation of $\omega_c$. This argument requires that $\Phi_t$ be a continuous bijection, so that $\Phi_t(D)$ remains a well-defined surface with ${\partial(\Phi_t(D)) = \Phi_t(\gamma)}$. Barham and Burby \cite{barham_diagnosing_2025} turn this into a numerical diagnostic for symplecticity: sample $N_s$ points around $\gamma$, evolve them one timestep, and check whether the discretized loop integral is conserved.

\subsection{PIC Hamiltonian, splitting, and the $C^1$ requirement}

The phase space of a discrete Vlasov--Maxwell system can be defined by the following canonical coordinates, corresponding to $N_\text{part}$ particles coupled to finite element field DOFs:
\begin{eqn}
\w{z}=(\w{a},\w{x}_1,\dots,\w{x}_{N_\text{part}};~\w{e},\w{p}_1,\dots,\w{p}_{N_\text{part}}),
\end{eqn}
where ${(\w{a},\w{e})}$ are coefficients for the finite element vector potential and electric field, as in Eq.~(\ref{discHamilSystem}), and $(\w{x}_s,\w{p}_s)$ are the position and canonical momentum of particle $s$. The Poisson bracket and the Hamiltonian, approximated with $\mQ_1$ in place of $\mM_1^{-1}$, extend the vacuum case of Eq.~(\ref{discHamilSystem}), namely
\begin{eqn}
\{F,G\}_\text{PIC}&=\frac{1}{\epsilon_0}\left(\pd{F}{\w{e}}\cdot\pd{G}{\w{a}}-\pd{G}{\w{e}}\cdot\pd{F}{\w{a}}\right)+\sum\limits_{s=1}^{N_\text{part}}\left(\pd{F}{\w{x}_s}\cdot\pd{G}{\w{p}_s}-\pd{G}{\w{x}_s}\cdot\pd{F}{\w{p}_s}\right)\\
\tilde{H}_\text{PIC}&=\underbrace{\frac{\epsilon_0}{2}\w{e}^T\mQ_1\w{e}\rule{0pt}{16pt}}_{\tilde{H}_\w{e}}+\underbrace{\frac{1}{2\mu_0}\w{a}^T\mC^T\tilde{\mM}_2\mC\,\w{a}\rule{0pt}{16pt}}_{\tilde{H}_\w{a}}+\underbrace{\sum\limits_{s=1}^{N_\text{part}}\frac{1}{2m_s}\abs{\w{p}_s-q_s\w{a}\cdot\gv{\Lambda}^1(\w{x}_s)}^2\rule{0pt}{16pt}}_{\tilde{H}_\text{Kin}=\tilde{H}_\text{Kin}^x+\tilde{H}_\text{Kin}^y+\tilde{H}_\text{Kin}^z},
\label{PICHamiltonian}
\end{eqn}
where $\gv{\Lambda}^1(\w{x})$ is the vector of 1-form basis functions evaluated at ${\w{x}\in\mR^3}$, so that ${\w{A}(\w{x})=\w{a}\cdot\gv{\Lambda}^1(\w{x})}$ is the interpolated vector potential. 

This Hamiltonian admits a natural five-way splitting \cite{he_hamiltonian_2015,glasser_gauge-compatible_2022} into $\tilde{H}_\w{e}$, $\tilde{H}_\w{a}$, and $\tilde{H}_\text{Kin}^\alpha$ for ${\alpha\in\{x,y,z\}}$, each of which is exactly integrable, with flow maps that can be evaluated in closed form and applied numerically to machine precision. Their respective EOMs (omitting static DOFs in each subsystem) are:
\begin{eqn}
\tilde{H}_\w{e}:&\quad \dot{\w{a}} = -\mQ_1\w{e}\\
\tilde{H}_\w{a}:&\quad \dot{\w{e}} = c^2\mC^T\tilde{\mM}_2\mC\,\w{a}\\
\tilde{H}_\text{Kin}^\alpha:&\quad\begin{cases}
\dot{x}^\alpha_s = \dfrac{1}{m_s}\bigl(p_{s\alpha}-q_s A(\w{x}_s)_\alpha\bigr)\\[8pt]
\dot{p}_{s\mu} = q_s\,\dot{x}^\alpha_s\;\partial_{x^\mu}A(\w{x}_s)_\alpha~~~~(\forall~~\mu\in\{x,y,z\})\\[8pt]
\dot{\w{e}} = -\frac{1}{\epsilon_0}\displaystyle\sum_{s=1}^{N_\text{part}} q_s\,\dot{x}^\alpha_s\;\gv{\Lambda}^1(\w{x}_s)_{\alpha}.
\end{cases}
\label{allEOM}
\end{eqn}
The subsystems $\tilde{H}_\w{e}$ and $\tilde{H}_\w{a}$ are linear and pose no smoothness difficulties. The kinetic subsystem $\tilde{H}_\text{Kin}^\alpha$ couples particles to the field through the interpolated vector potential, and the momentum equation involves $\partial_{x^\mu}A(\w{x}_s)_\alpha$---the spatial gradient of the 1-form interpolant evaluated at the particle position. This term is the origin of the $C^1$ smoothness requirement. Note that no summation convention is used in Eq.~(\ref{allEOM}); within each sub-Hamiltonian $\tilde{H}_\text{Kin}^\alpha$, the index $\alpha\in{x,y,z}$ is fixed, and $\mu$ takes each value ${x,y,z}$ in turn.

During each kinetic sub-flow (i.e. in any of the three ${\tilde{H}_\text{Kin}^\alpha}$ substeps), the vector potential DOFs $\w{a}$ are frozen. The particles evolve at fixed field configuration, with the force in Eq.~(\ref{allEOM}) determined by $\partial_{x^\mu}A(\w{x}_s)_\alpha$ at the instantaneous particle position. If the 1-form basis functions $\gv{\Lambda}^1$ are only $C^0$ across cell boundaries, then $\partial_{x^\mu}A$ has jump discontinuities at cell faces. The kinetic sub-flow therefore inherits these discontinuities: two particles at positions infinitesimally apart but on opposite sides of a cell face experience different forces for the duration of the sub-step, yielding a spatially discontinuous map. Phase-space loops that straddle a cell face are torn, and the topological premise of the Stokes' theorem argument for Eq.~(\ref{PoincareInvariant}) fails. Barham and Burby confirm this mechanism analytically for a model problem and numerically for an electrostatic PIC simulation; the argument generalizes to the full electromagnetic case through the kinetic sub-Hamiltonians ${\tilde{H}_\text{Kin}^\alpha}$ above, whose discontinuity mechanism is identical.

The EOMs of Eq.~(\ref{allEOM}) make the smoothness requirement explicit: the Hamiltonian vector field involves $\partial_{x^\mu}A(\w{x})_\alpha$, i.e.\ the first spatial derivatives of the 1-form basis functions. For this vector field to be spatially $C^0$---and hence for the splitting map to be continuous---the 1-form basis must be at least $C^1$.

For degree-$p$ B-spline discretizations on cubical meshes \cite{buffa_isogeometric_2010,buffa_isogeometric_2011,perse_geometric_2021}, the basis functions are $C^{p-1}$ uniformly across all form degrees, so the $C^1$ requirement on 1-forms is met for B-splines of degree ${p\geq 2}$. On simplicial meshes, smooth de~Rham complexes (Fig.~\ref{deRhamDiagram}) with sufficient regularity have been constructed \cite{christiansen_nodal_2018,fu_alfeld_2020} but not, to our knowledge, been applied to PIC. Whitney forms (tangential continuity only; not $C^0$ pointwise), $C^0$ Lagrange elements, and $C^0$ spectral elements all fail this requirement.

B-spline GYMs with $p\geq 2$ therefore satisfy both structural requirements identified in this paper: the sparse approximate inverse mass matrices of Section~\ref{SparseMassMatSect} restore the scalability of the field solve without affecting symplecticity (Theorem~\ref{thm:symplecticity}), while the $C^{p-1}$ smoothness of the basis ensures that the Hamiltonian splitting flow is spatially continuous and preserves the Poincar\'e integral invariant over particle trajectories.

For PIC applications, SPAI-OP (Section~\ref{SPAIOPsubsect}) further supports targeted dispersion control: probing vectors chosen from the modes most responsible for numerical Cherenkov radiation~\cite{godfrey_cherenkov_2014} can concentrate accuracy where it is most needed for beam-driven instability suppression, without sacrificing the $C^{p-1}$ smoothness guaranteed by high-order B-splines.

The gauge-compatibility structure of Section~\ref{GYMdefnSubsect} extends naturally to the PIC setting. Whereas the vacuum gauge transformation $\w{a}\mapsto\w{a}+\mG\w{s}$ leaves $\tilde{H}\w{a}$ invariant through the identity ${\mC\mG=\mZero}$, the particle coupling in $\tilde{H}\text{Kin}$ requires a compensating shift of the canonical momenta to preserve the Hamiltonian. The extended gauge transformation is
\begin{eqn}
\left(\begin{matrix}\w{a}\\\w{p}_s\end{matrix}\right)\mapsto\left(\begin{matrix}\hspace{-48pt}\w{a}+\mG\w{s}\\\w{p}_s+q_s\mG\w{s}\cdot\gv{\Lambda}^1(\w{x}_s)\end{matrix}\right)
\end{eqn}
$\forall$ ${1\leq s\leq N_\text{part}}$. As shown in \cite{glasser_gauge-compatible_2022}, this transformation is generated by the Gauss' law momentum map
\begin{eqn}
\mu=\mG^T\w{e}+\sum\limits_s\frac{q_s}{\epsilon_0}\gv{\Lambda}^0(\w{x}_s).
\end{eqn}
Physically, the momentum map $\mu$ is the discrete version of Gauss' law, ${\nabla\cdot\w{E} - \rho/\epsilon_0}$, where $\mG^T\w{e}$ discretizes the divergence of $\w{E}$ and $\sum_s q_s \gv{\Lambda}^0(\w{x}_s)/\epsilon_0$ discretizes the charge density. $\mu$ is exactly preserved by the flow of each sub-Hamiltonian---${\{\mu,\tilde{H}_i\}_\text{PIC}=0}$ $\forall$ ${\tilde{H}_i\in\{\tilde{H}_\w{e},\tilde{H}_\w{a},\tilde{H}_\text{Kin}^\alpha\}}$---and therefore by the composed GYM flow, guaranteeing exact discrete charge conservation at every time step.

\section{Numerical Results\label{NumResult}}

\subsection{Error scaling with mesh refinement\label{errorScaling}}

\begin{figure*}[t!]
\centering
\includegraphics[width=\textwidth]{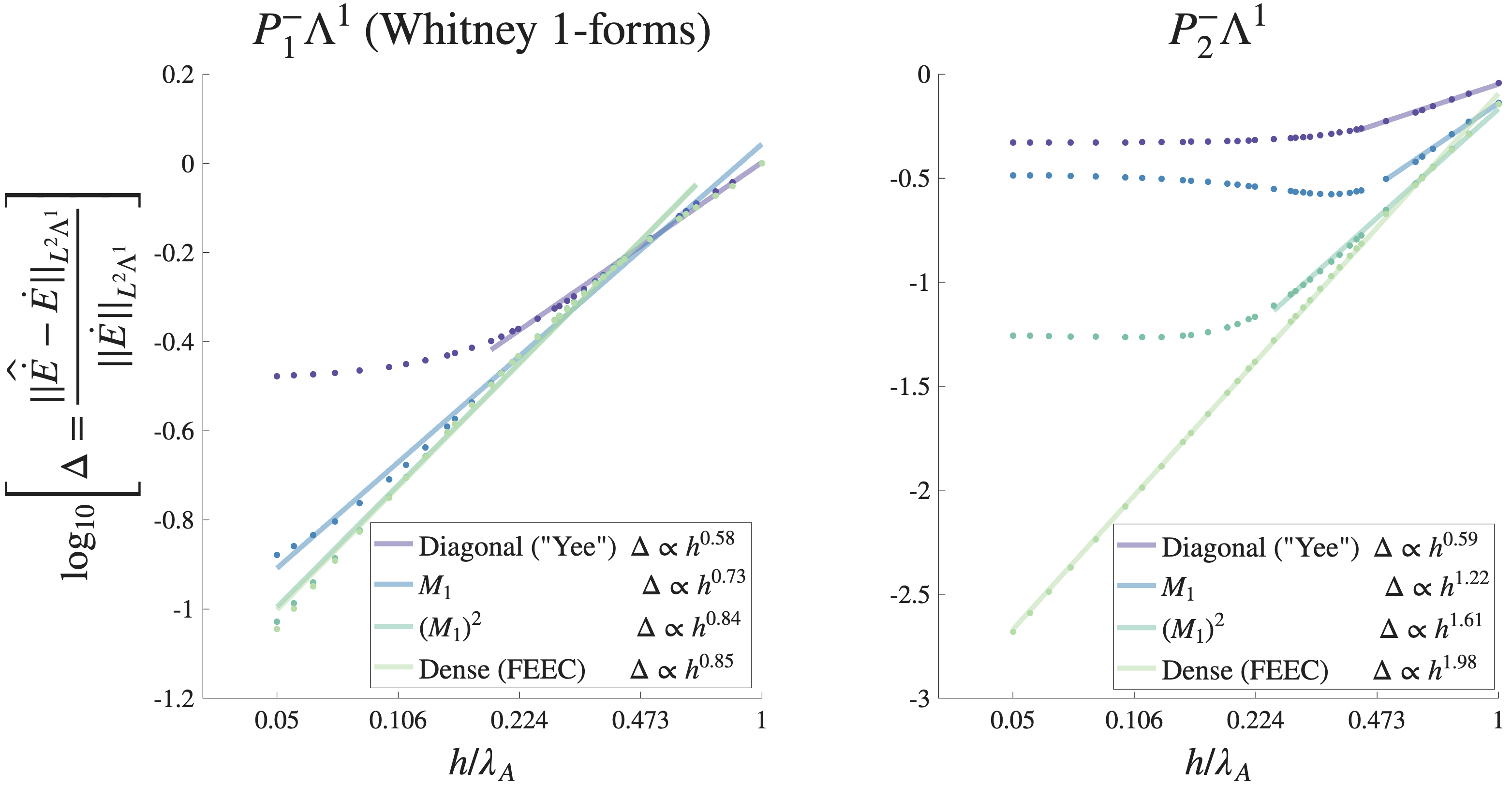}
\caption{These plots depict the error in various finite element approximations of the curl-of-curl operator. The relative ${L^2\Lambda^1}$ error in the discrete approximation ${\smash{\hat{\dot{\w{E}}}}}$ of ${\dot{\w{E}}=c^2\nabla\times\nabla\times\w{A}}$ is plotted against the cell size $h$ (measured per signal wavelength $\lambda_\w{A}$) of Delaunay triangulations of a 2D periodic domain. Four sparsity patterns for ${\mQ_1\approx\mM_1^{-1}}$ are compared---diagonal (``Yee''), $\mM_1$, ${(\mM_1)^2}$, and dense (exact FEEC)---for two finite element families: first-order Whitney 1-forms ${\mmP_1^-\Lambda^1}$ (left) and second-order trimmed polynomial 1-forms ${\mmP_2^-\Lambda^1}$ (right). Highlighted linear overlays indicate the convergent regime for each data series, with power-law scalings fit therein. The $\mM_1$-sparsity pattern achieves much of the improved scaling that exact FEEC provides for Whitney forms, while ${(\mM_1)^2}$-sparsity captures much of the second-order accuracy scaling. The error saturation at small ${h/\lambda_\w{A}}$ reflects Strang consistency error \cite{strang_variational_1972}, which is independent of mesh refinement, arising inevitably from ``localizing" (or sparsifying) the exact inverse $\mM_1^{-1}$. Wider stencils---$\mM_1$, $(\mM_1)^2$---reduce this consistency error by incorporating progressively more of the non-local inverse.}
\label{CurlOfCurlFig}
\end{figure*}

To measure the accuracy of the approximation ${\mQ_1\approx\mM_1^{-1}}$ in a manner relevant to a GYM, we first examine the curl-of-curl operator ${\nabla\times\nabla\times}$ (discretized by ${\mM_1^{-1}\mC^T\mM_2\mC}$ in the FEEC setting). We generate 2D simplicial Delaunay triangulations ${\mT_h}$ of the flat torus ${\abs{\mT_h}=[0,L_x]\times[0,L_y]\subset\mR^2}$ with periodic boundary conditions. On $\mT_h$, we consider a sinusoidal vector potential ${\w{A}=\sin(k_ny)\md x}$ for ${k_n=2\pi n/L_y}$ and ${n\in\mN}$, and canonically project it onto a discrete FEEC basis \cite{arnold_finite_2006}, ${\w{A}\mapsto\w{a}\cdot\gv{\Lambda}^1\in\Lambda^1(\mT_h)}$. We measure the relative ${L^2\Lambda^1}$ error
\begin{eqn}
\Delta \equiv \norm{\hat{\dot{\w{E}}}-\dot{\w{E}}}_{L^2\Lambda^1}\Big/\norm{\dot{\w{E}}}_{L^2\Lambda^1}
\label{relErrorDef}
\end{eqn}
between the continuum curl-of-curl ${\dot{\w{E}}/c^2=\nabla\times\nabla\times\w{A}=k_n^2\sin(k_ny)\md x}$ and its discrete approximation ${\hat{\dot{\w{E}}}/c^2=\mQ_1\mC^T\mM_2\mC\w{a}\cdot\gv{\Lambda}^1}$, for various sparsity patterns ${S(\mQ_1)}$. The results are shown in Fig.~\ref{CurlOfCurlFig}.

The figure compares two FEEC bases---the first-order Whitney 1-forms ${\mmP_1^-\Lambda^1(\mT_h)}$ and the second-order trimmed polynomial 1-forms ${\mmP_2^-\Lambda^1(\mT_h)}$---across four sparsity patterns for $\mQ_1$: diagonal (Yee's implicit pattern), $\mM_1$ sparsity, ${(\mM_1)^2}$ sparsity, and dense (the exact inverse). Each sparse $\mQ_1$ is computed via the column-wise Frobenius-optimal SPAI of Eq.~(\ref{NormalEqSolve}), followed by symmetrization ${\mQ_1\leftarrow\frac{1}{2}(\mQ_1+\mQ_1^T)}$.

For Whitney forms (left panel of Fig.~\ref{CurlOfCurlFig}), the symmetrized Frobenius-optimal SPAI at $\mM_1$-sparsity achieves an error scaling of ${h^{0.73}}$, compared to ${h^{0.85}}$ for the exact FEEC operator and ${h^{0.58}}$ for the diagonal (Yee) pattern. Beyond the improvement in scaling exponent, the $\mM_1$-sparsity pattern critically extends the convergent regime. The diagonal approximation saturates at ${h/\lambda_\w{A}\lesssim 0.2}$ (where ${\lambda_\w{A}=L/n}$ denotes signal wavelength); beyond this point, further mesh refinement yields no improvement, because the error is dominated by the \emph{consistency error} of the sparse approximation $\mQ_1\approx\mM_1^{-1}$ rather than the finite element interpolation error. This saturation is a manifestation of replacing the exact bilinear form (which involves $\mM_1^{-1}$) with an approximate one (involving $\mQ_1$), introducing an error term in the first Strang Lemma that is independent of mesh size $h$ \cite{strang_variational_1972,strang_analysis_1973}. The exact inverse $\mM_1^{-1}$ is a non-local operator---its Green's function decays but does not vanish with distance---so a fixed-stencil sparse approximation captures only a fixed fraction of its action, regardless of resolution. The $\mM_1$- and ${(\mM_1)^2}$-sparsity patterns delay this saturation by incorporating progressively more non-local coupling. For second-order elements (right panel), ${(\mM_1)^2}$-sparsity achieves ${h^{1.61}}$ scaling versus ${h^{1.98}}$ for exact FEEC, likewise recovering much of the higher-order accuracy scaling.

The additional computational cost of moving from diagonal to $\mM_1$-sparsity is modest: the communicating node pairs remain unchanged, with roughly 5 times more data per message in our 2D example. Since communication volume scales with the number of boundary cells, the overhead is small relative to the accuracy gained.

\subsection{SPAI-OP eigenvalue accuracy on structured grids\label{SPAIOPresults}}

We next evaluate the SPAI-OP formulation of Section~\ref{SPAIOPsubsect} on a structured periodic grid. We construct a ${20\times 80}$ flat torus ${\abs{\mT_h}=[0,5]\times[0,20]\subset\mR^2}$ with square cells (${h_x=h_y=0.25}$) and periodic boundary conditions. The 1-form mass matrix $\mM_1$ corresponds to the Whitney $Q_1^-\Lambda^1$ formula on rectangles ($1/3$ on the diagonal, $1/6$ off-diagonal for each face pair sharing a cell), and we continue to denote ${\mP_2=\mC^T\mM_2\mC}$. Three methods are compared at each sparsity pattern ${S(\mQ_1)=S(\mM_1^k)}$ for ${k\in\{0,1,2\}}$:
\begin{enumerate}
\item \textbf{SPAI}: the column-wise Frobenius-optimal solution of Eq.~(\ref{NormalEqSolve}), followed by post-symmetrization.  On the structured periodic grid used here the column-wise SPAI is already self-symmetric, so this coincides with the Sylvester $\lambda{=}0$ solution of Eq.~(\ref{SPAIOPSylvester}).
\item \textbf{SPAI-OP at} ${\lambda^*(\alpha=10\%)}$: the operator-probed formulation, Eq.~(\ref{SPAIOPSylvester}) with ${\hat{\mmT}=\mOne+\lambda\mmW\mmW^T}$, probing the ${n_\text{probe}=100}$ lowest non-kernel eigenmodes of ${\mK=\mM_1^{-1}\mP_2}$, i.e., ${\mmW=\mP_2\mV}$ with ${\mV=[\w{v}_1,\dots,\w{v}_{100}]}$ as defined in Section~\ref{SPAIOPsubsect}.  The weight ${\lambda^*(\alpha)}$ is the tolerance-optimal value defined in Eq.~(\ref{kneeCriterion}), explored in Section~\ref{lambdaScan}.
\item \textbf{SPAI-OP at} ${\lambda^*(\alpha=25\%)}$: the same formulation at the paper-default tolerance.
\end{enumerate}
Accuracy is reported as the residual eigenvalue error ${\varepsilon_k=\norm{\tilde{\mK}\w{v}_k-\lambda_k\w{v}_k}/\abs{\lambda_k}}$ defined in Eq.~(\ref{eigError})  per mode, which measures how faithfully the approximate curl-of-curl operator preserves each exact eigenmode.

SPAI-OP concentrates accuracy on probed modes. Table~\ref{SPAIOPeigTable} reports the mode-averaged residual eigenvalue error of Eq.~(\ref{eigError}) for three methods---the symmetric column-wise SPAI of Eq.~(\ref{NormalEqSolve}); SPAI-OP at the light-tolerance ${\lambda^*(\alpha=10\%)}$; and SPAI-OP at the paper-default ${\lambda^*(\alpha=25\%)}$---with means separately computed over the full non-kernel spectrum, the $n_\text{probe}=100$ probed modes, and the remaining $1499$ unprobed modes.  The paper-default row is the principal result: across all three sparsity patterns, SPAI-OP at $\lambda^*(25\%)$ reduces the mean error on the probed modes by $2.2\times$ at $S(\mM_1)$, $2.2\times$ at $S(\mM_1^2)$, and $7.8\times$ at $S(\mOne)$ (where the tolerance never binds and $\lambda^*$ sits at the scan-range saturation point), at the cost of a uniform $\sim 1.3\times$ increase in the mean error of the unprobed modes.  The light-tolerance value $\lambda^*(10\%)$ offers an intermediate operating point with $1.4$--$2.0\times$ probed-mode gain and a smaller unprobed-error penalty ($\sim 1.1\times$).

\begin{table}[h!]
\centering
\caption{Mode-averaged residual eigenvalue error $\bar\varepsilon=\mathrm{mean}_k\norm{\tilde{\mK}\w{v}_k-\lambda_k\w{v}_k}/|\lambda_k|$ on the ${20\times80}$ periodic structured grid for SPAI and SPAI-OP evaluated at the tolerance-optimal $\lambda^*$'s of Table~\ref{kneeTable} (probe set: $100$ lowest non-kernel modes).  Means are reported over all non-kernel modes, the $100$ probed modes, and the $1499$ unprobed modes.  Bold rows correspond to the paper default $\alpha=25\%$.}
\label{SPAIOPeigTable}
\begin{tabular}{llccc}
\hline
Pattern & Method & $\bar\varepsilon$ (all) & $\bar\varepsilon$ (probed) & $\bar\varepsilon$ (unprobed) \\
\hline
$S(\mOne)$    & SPAI                                    & $0.3394$ & $0.2714$ & $0.3439$ \\
              & SPAI-OP $\lambda^*(10\%)$               & $0.3731$ & $0.1383$ & $0.3888$ \\
              & \textbf{SPAI-OP $\lambda^*(25\%)$}      & $\mathbf{0.4185}$ & $\mathbf{0.0347}$ & $\mathbf{0.4441}$ \\
\hline
$S(\mM_1)$    & SPAI                                    & $0.0955$ & $0.0493$ & $0.0986$ \\
              & SPAI-OP $\lambda^*(10\%)$               & $0.1045$ & $0.0295$ & $0.1095$ \\
              & \textbf{SPAI-OP $\lambda^*(25\%)$}      & $\mathbf{0.1192}$ & $\mathbf{0.0224}$ & $\mathbf{0.1257}$ \\
\hline
$S(\mM_1^2)$  & SPAI                                    & $0.0259$ & $0.0146$ & $0.0266$ \\
              & SPAI-OP $\lambda^*(10\%)$               & $0.0284$ & $0.0106$ & $0.0296$ \\
              & \textbf{SPAI-OP $\lambda^*(25\%)$}      & $\mathbf{0.0323}$ & $\mathbf{0.0065}$ & $\mathbf{0.0340}$ \\
\hline
\end{tabular}
\end{table}

Per mode, SPAI-OP concentrates its error reduction sharply within the probed subset (${k\leq 100}$), while the 1499 unprobed modes absorb the redistributed cost summarized in Table~\ref{SPAIOPeigTable}.

\subsection{Probing weight sensitivity and the tolerance-optimal $\lambda^*$\label{lambdaScan}}

\begin{figure*}[t!]
\centering
\includegraphics[width=\textwidth]{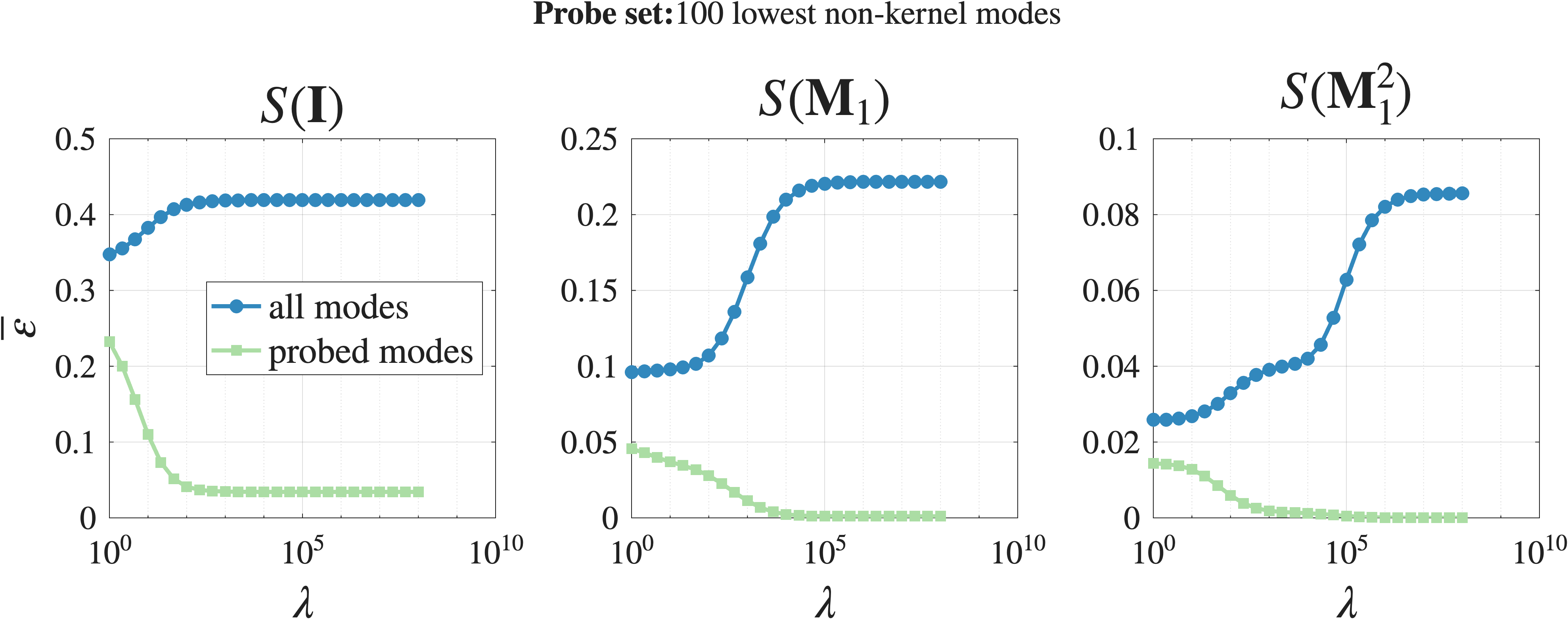}
\caption{SPAI-OP probing-weight sensitivity on the ${20\times80}$ structured grid, bulk-concentration probe set ($n_\text{probe}=100$ lowest modes).  Blue circles: mean eigenvalue error over all non-kernel modes.  Green squares: mean error over the probed subset.  Both traces share a linear $\bar\varepsilon$ axis; $\lambda$ is on a log axis.  As $\lambda$ increases, probed-mode error decreases monotonically while all-mode error rises monotonically toward a saturation plateau---an accuracy-vs-coverage tradeoff.}
\label{lambdaScanFig}
\end{figure*}

The probing weight $\lambda$ is the single tunable parameter of SPAI-OP, interpolating between pure symmetric SPAI (${\lambda=0}$) and strongly mode-concentrated accuracy (${\lambda\gg1}$).  A principled choice of $\lambda$ is therefore the main practical question for a user of the method.  The lambda scan of Fig.~\ref{lambdaScanFig} records the mean eigenvalue error over (i) all non-kernel modes, $\bar\varepsilon_\text{all}(\lambda)$, and (ii) the probed subset, $\bar\varepsilon_\text{probed}(\lambda)$, for the three sparsity patterns ${S(\mM_1^k)}$ and two representative probe configurations:
\begin{itemize}
\item \textbf{bulk concentration}: the $n_\text{probe}=100$ lowest non-kernel modes, suitable when the practitioner wishes to protect the full resolvable wave band (e.g.\ broadband simulations in which multiple long-wavelength modes must propagate accurately);
\item \textbf{targeted single-frequency}: $n_\text{probe}=5$ modes drawn from mid-spectrum (indices $797{-}801$ of $1599$, eigenvalues $137.9{-}138.3$), suitable when the simulation concerns one physical wavelength (e.g.\ a narrow-band wave packet or a beam-driven instability).
\end{itemize}
In both cases the eigenvalue error is computed on each \emph{physical} eigenmode via the residual metric
\begin{eqn}
\varepsilon_k\equiv\frac{\norm{\tilde{\mK}\w{v}_k-\lambda_k\w{v}_k}}{|\lambda_k|},
\label{eigError}
\end{eqn}
where $\w{v}_k$ is the $k^\text{th}$ exact non-kernel eigenvector of $\mK=\mM_1^{-1}\mP_2$ with eigenvalue $\lambda_k$.  The residual metric is mode-faithful: it measures how well $\tilde{\mK}=\mQ_1\mP_2$ maps each exact eigenvector to $\lambda_k\w{v}_k$, regardless of how the approximate spectrum reorders.

\begin{figure*}[t!]
\centering
\includegraphics[width=\textwidth]{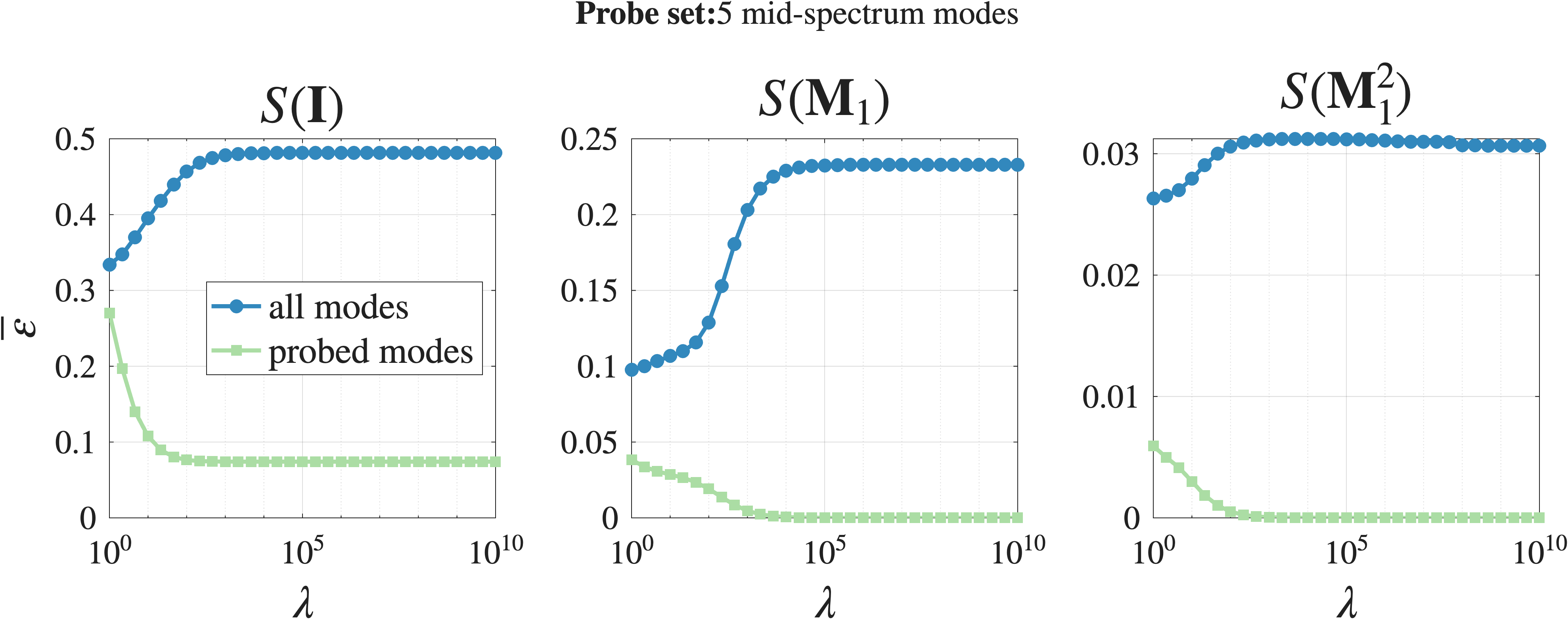}
\caption{Same as Fig.~\ref{lambdaScanFig} but for the targeted single-frequency probe set ($n_\text{probe}=5$ mid-spectrum modes).  Note that for ${S(\mM_1^2)}$ sparsity the probed-mode error falls to machine precision while the all-mode error increases only slightly, and that the ${S(\mM_1)}$ and ${S(\mM_1^2)}$ curves drive the probed error essentially to zero, reflecting the abundance of $\mQ_1$ degrees of freedom relative to the small number of probing constraints.}
\label{lambdaScanMid5Fig}
\end{figure*}

\paragraph{Tolerance-optimal $\lambda^*$ criterion}  A principled choice of $\lambda$ follows from bounding the all-mode degradation.  We define the tolerance-optimal value
\begin{eqn}
\lambda^*(\alpha) \equiv \max\,\bigl\{\,\lambda \,:\, \bar\varepsilon_\text{all}(\lambda) \leq (1+\alpha)\,\bar\varepsilon_\text{all}^\text{SPAI}\,\bigr\},
\label{kneeCriterion}
\end{eqn}
where $\alpha$ is a user-chosen tolerance (paper default $\alpha=25\%$).  Because $\bar\varepsilon_\text{all}$ is already a per-mode average over the full non-kernel spectrum, bounding its relative increase directly controls the average per-mode degradation across \emph{all} modes, automatically weighting each mode equally regardless of the probed-subset size.  The tolerance $\alpha$ is thus interpretable as a budget on per-mode SPAI-baseline degradation, independent of how many (or few) modes are probed.

Table~\ref{kneeTable} reports $\lambda^*(\alpha)$ for $\alpha\in\{10\%,25\%,50\%\}$ across the six combinations of probe set and sparsity pattern.  The table also reports the resulting probed-mode error $\bar\varepsilon_\text{probed}(\lambda^*)$ and the \emph{gain factor} ${\bar\varepsilon_\text{probed}^\text{SPAI}/\bar\varepsilon_\text{probed}(\lambda^*)}$, which quantifies how much better the probed modes are resolved relative to pure SPAI.

\begin{table*}[h!]
\centering
\caption{Tolerance-optimal $\lambda^*(\alpha)$ for SPAI-OP on the ${20\times80}$ structured periodic grid, defined by Eq.~(\ref{kneeCriterion}).  ``Gain'' is the ratio $\bar\varepsilon_\text{probed}^\text{SPAI}/\bar\varepsilon_\text{probed}(\lambda^*)$; ``sat.''\ denotes that the all-mode tolerance is not binding within the scan range, in which case the probed-mode error has effectively saturated and further $\lambda$ buys no additional benefit.  Paper default is $\alpha=25\%$ (bold).}
\label{kneeTable}
\begin{tabular}{llccccc}
\hline
Probe set & Pattern & $\alpha$ & $\lambda^*$ & $\bar\varepsilon_\text{all}(\lambda^*)$ & $\bar\varepsilon_\text{probed}(\lambda^*)$ & gain \\
\hline
\multirow{9}{*}{\parbox{1.9cm}{bulk, \\ $n_\text{probe}=100$}}
 & $S(\mOne)$   & $10\%$ & $6.2{\cdot}10^{0}$ & $0.373$ & $0.137$ & $2.0\times$ \\
 & $S(\mOne)$   & $25\%$ & $\geq 10^{8}$ (sat.) & $0.419$ & $0.034$ & $7.9\times$ \\
 & $S(\mOne)$   & $50\%$ & $\geq 10^{8}$ (sat.) & $0.419$ & $0.034$ & $7.9\times$ \\
 \cline{2-7}
 & $S(\mM_1)$   & $10\%$ & $7.5{\cdot}10^{1}$ & $0.105$ & $0.029$ & $1.7\times$ \\
 & $S(\mM_1)$   & \textbf{25\%} & $\mathbf{2.2{\cdot}10^{2}}$ & $\mathbf{0.119}$ & $\mathbf{0.022}$ & $\mathbf{2.2\times}$ \\
 & $S(\mM_1)$   & $50\%$ & $5.9{\cdot}10^{2}$ & $0.143$ & $0.015$ & $3.3\times$ \\
 \cline{2-7}
 & $S(\mM_1^2)$ & $10\%$ & $2.5{\cdot}10^{1}$ & $0.028$ & $0.011$ & $1.4\times$ \\
 & $S(\mM_1^2)$ & $25\%$ & $8.5{\cdot}10^{1}$ & $0.032$ & $0.006$ & $2.3\times$ \\
 & $S(\mM_1^2)$ & $50\%$ & $8.5{\cdot}10^{2}$ & $0.039$ & $0.002$ & $7.3\times$ \\
\hline
\multirow{9}{*}{\parbox{1.9cm}{targeted, \\ $n_\text{probe}=5$}}
 & $S(\mOne)$   & $10\%$ & $5.1{\cdot}10^{0}$ & $0.373$ & $0.135$ & $3.1\times$ \\
 & $S(\mOne)$   & $25\%$ & $2.7{\cdot}10^{1}$ & $0.424$ & $0.087$ & $4.8\times$ \\
 & $S(\mOne)$   & $50\%$ & $\geq 10^{10}$ (sat.) & $0.482$ & $0.074$ & $5.6\times$ \\
 \cline{2-7}
 & $S(\mM_1)$   & $10\%$ & $6.8{\cdot}10^{0}$ & $0.105$ & $0.030$ & $2.1\times$ \\
 & $S(\mM_1)$   & \textbf{25\%} & $\mathbf{5.8{\cdot}10^{1}}$ & $\mathbf{0.119}$ & $\mathbf{0.022}$ & $\mathbf{2.8\times}$ \\
 & $S(\mM_1)$   & $50\%$ & $1.6{\cdot}10^{2}$ & $0.143$ & $0.016$ & $4.0\times$ \\
 \cline{2-7}
 & $S(\mM_1^2)$ & $10\%$ & $1.4{\cdot}10^{1}$ & $0.028$ & $0.002$ & $7.4\times$ \\
 & $S(\mM_1^2)$ & $25\%$ & $\geq 10^{10}$ (sat.) & $0.031$ & $\approx 0$ & $>10^{4}\times$ \\
 & $S(\mM_1^2)$ & $50\%$ & $\geq 10^{10}$ (sat.) & $0.031$ & $\approx 0$ & $>10^{4}\times$ \\
\hline
\end{tabular}
\end{table*}

Two regimes emerge clearly from Table~\ref{kneeTable}. In the \emph{bulk} regime, the binding $\lambda^*(25\%)$ values are of order $100$ across sparsity patterns (with $S(\mOne)$ a special case in which the tolerance never binds, discussed below), delivering $2.2$--$7.9\times$ probed-mode gain over pure SPAI.  The $S(\mOne)$ (Yee-like) case is special: the all-mode error saturates below the $25\%$ threshold, so the tolerance never binds; a practitioner can push $\lambda$ past the saturation point ($\lambda\gtrsim 100$) and extract the full $\sim 8\times$ benefit at no tolerance cost.  For $S(\mM_1)$ and $S(\mM_1^2)$ the tolerance is binding, and $\alpha$ directly controls the position on a continuous accuracy/degradation frontier.

In the \emph{targeted} regime (few probes), the probing constraint is cheap to satisfy---only five equations per mode over thousands of degrees of freedom in $\mQ_1$---and the tolerance-optimal $\lambda^*$ values drop by two to three orders of magnitude.  The $S(\mM_1^2)$ case becomes qualitatively different: \emph{the all-mode tolerance never binds even at $\alpha=10\%$}, because the rich sparsity pattern can accommodate the five probing constraints to machine precision without meaningfully disturbing the remaining modes.  In this limit, SPAI-OP offers essentially free accuracy at the target frequency.  This is precisely the regime relevant to a single-frequency wave simulation---the signal practical application, examined quantitatively in Section~\ref{eigenmodeDispersion}.

\subsection{Stability and CFL verification\label{CFLverification}}

Proposition~\ref{prop:stability} predicts that the one-step Strang map is stable if and only if the timestep $\tau$ satisfies ${\tau\leq 2/\omega_{\max}}$, where $\omega_{\max}$ is the largest eigenfrequency of the approximate curl-of-curl operator ${\tilde{\mK}=\mQ_1\mP_2}$. We verify this for three choices of $\mQ_1$ approximation (diagonal, SPAI, SPAI-OP) on the ${20\times80}$ structured grid. For each method, we compute all eigenvalues of the one-step Jacobian $\mJ_{\Phi_\tau}$ at two timesteps: ${\tau=0.9\,\tau_{\text{CFL}}}$ (stable) and ${\tau=1.1\,\tau_{\text{CFL}}}$ (unstable), where ${\tau_{\text{CFL}}=2/\omega_{\max}}$.

\begin{figure*}[t!]
\centering
\includegraphics[width=\textwidth]{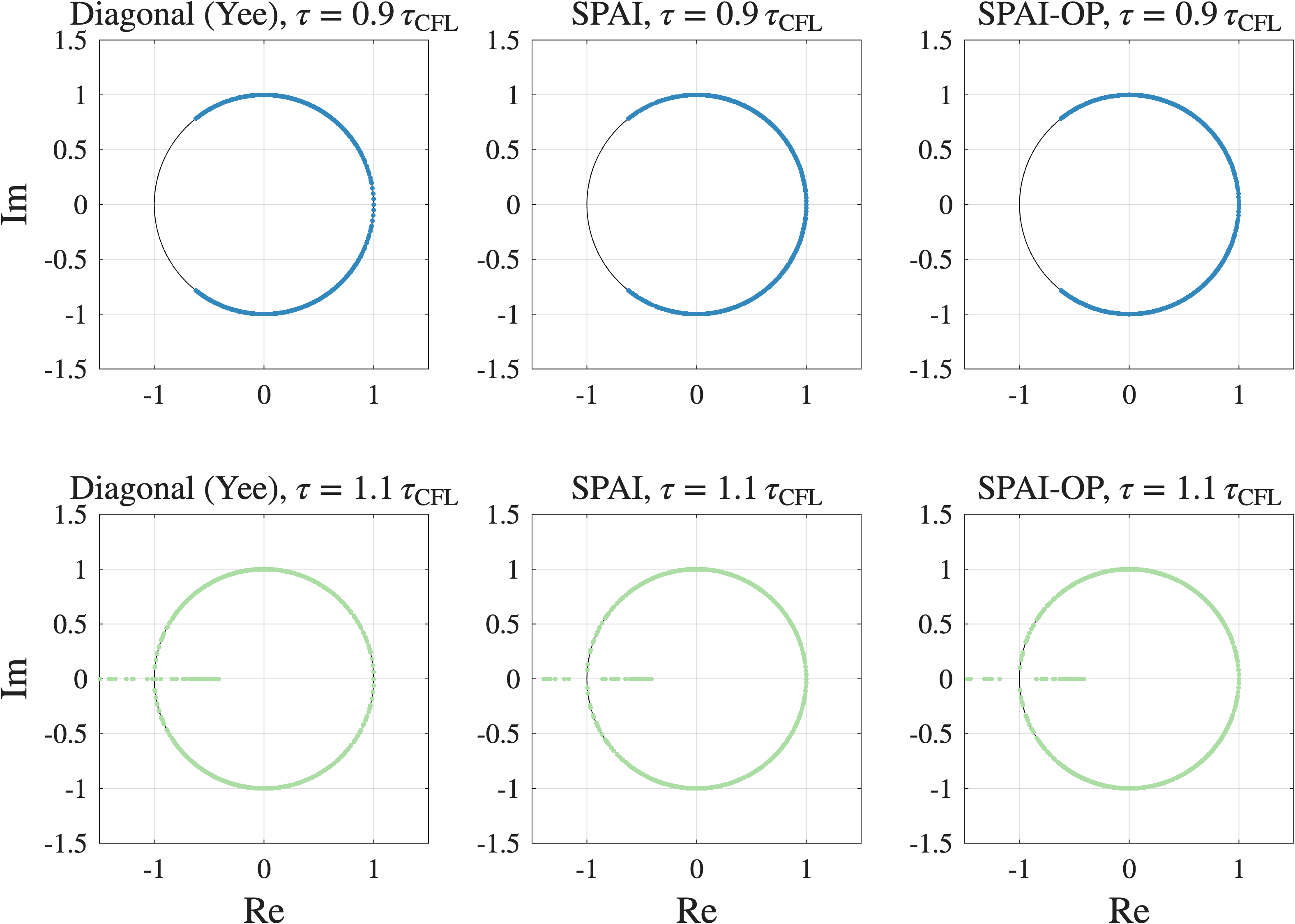}
\caption{Eigenvalues of the one-step Strang map $\Phi_\tau$ in the complex plane for three GYM configurations (diagonal, symmetric SPAI, SPAI-OP). Top row: ${\tau=0.9\,\tau_{\text{CFL}}}$ (stable); all eigenvalues lie on the unit circle. Bottom row: ${\tau=1.1\,\tau_{\text{CFL}}}$ (unstable); eigenvalue pairs leave the unit circle along the real axis. Each method has a different $\tau_{\text{CFL}}$ because different $\mQ_1$ produce different $\omega_{\max}$. The map is symplectic in both cases (eigenvalues come in reciprocal pairs); stability is a separate condition enforced by the CFL bound.}
\label{CFLfig}
\end{figure*}

Fig.~\ref{CFLfig} confirms the prediction: at $\tau<\tau_{\text{CFL}}$, all eigenvalues lie on the unit circle (top row), while at $\tau>\tau_{\text{CFL}}$, eigenvalue pairs depart along the real axis (bottom row). Crucially, the map remains symplectic in both cases---eigenvalues come in reciprocal pairs $(\lambda,1/\lambda)$, a hallmark of symplecticity---but stability requires the CFL condition. An indefinite $\mQ_1$ (not shown) produces eigenvalues off the unit circle at \emph{any} timestep, confirming the structural instability predicted by the analysis of Section~\ref{stabilitySect}.

\subsection{Time-domain energy conservation\label{timeDomainEnergy}}

To briefly examine symplecticity in a full time-domain simulation, we propagate the lowest non-kernel eigenmode of ${\mM_1^{-1}\mP_2}$ on a ${20\times20}$ periodic cubical mesh for 200 wave periods using the Strang splitting of Eq.~(\ref{strangSplit}).  The precise setup---including the eigenmode computation and the initial condition ${\w{e}(0)=\w{v}}$, ${\w{b}(0)=\w{0}}$---is described in Section~\ref{eigenmodeDispersion}; the present subsection reports the energy-conservation diagnostic from that same run.  We use first-order cubical Whitney forms (${Q_1^-\Lambda^p}$) with $S(\mM_1)$ sparsity for SPAI and SPAI-OP, and ${S(\mOne)}$ sparsity for the Yee-like baseline.

\begin{figure}[h!]
\centering
\includegraphics[width=\columnwidth]{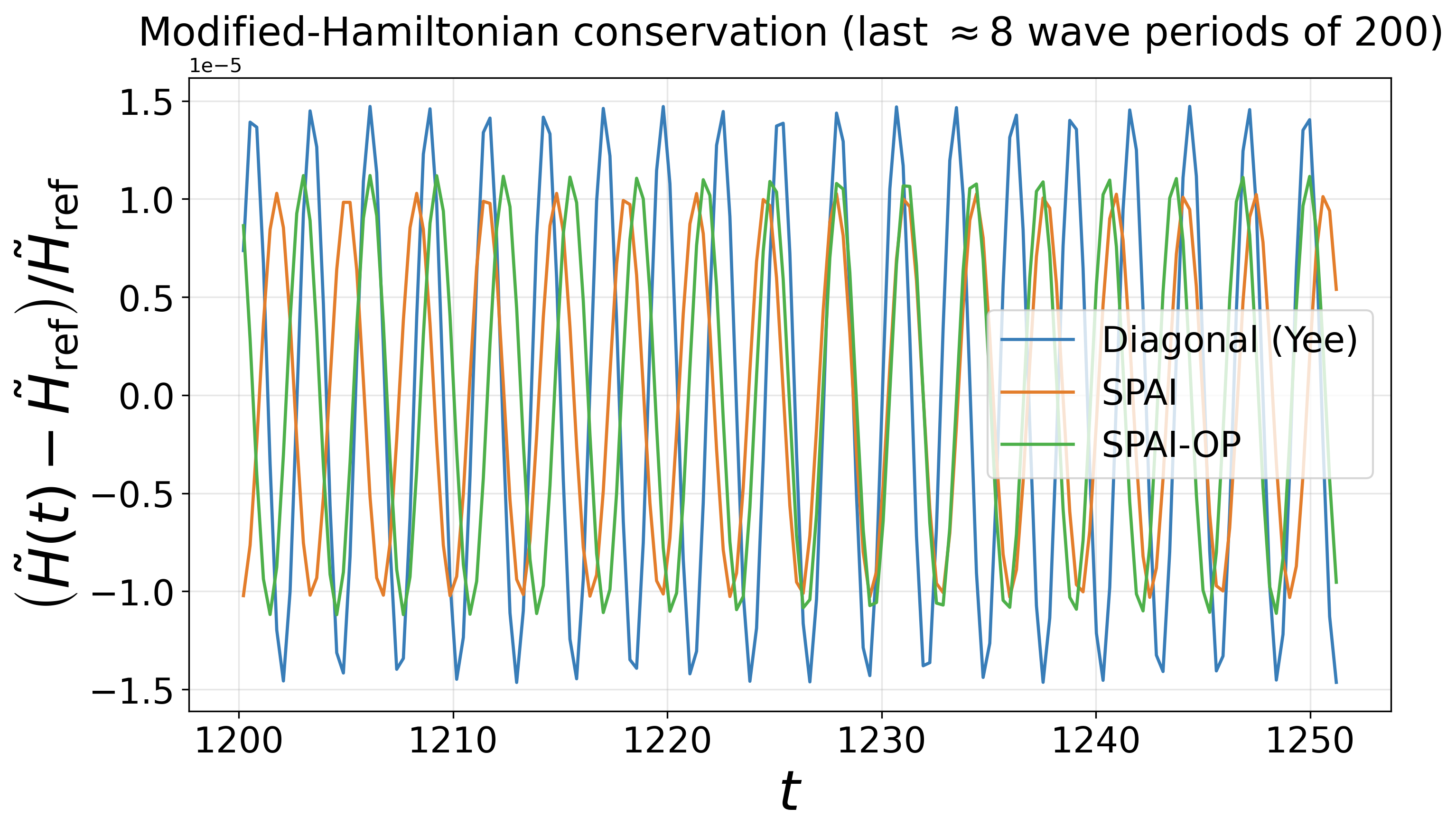}
\caption{Relative energy deviation over the final ${\sim}8$ wave periods of a 200-period simulation on a ${20\times20}$ cubical mesh with CFL safety factor $0.3$. All methods show bounded oscillation with no secular drift, consistent with the symplecticity of the Strang splitting proven in Theorem~\ref{thm:symplecticity}.}
\label{energyFig}
\end{figure}

The modified Hamiltonian ${\tilde{H}_\text{d}=\frac{1}{2}\w{e}^T\mQ_1\w{e}+\frac{1}{2}\w{b}^T\mM_2\w{b}}$ is recorded at each output step. Fig.~\ref{energyFig} shows the relative energy deviation ${(\tilde{H}(t)-\tilde{H}(0))/\tilde{H}(0)}$ over the full simulation.

All methods exhibit bounded energy oscillation with no secular drift, as expected from the symplecticity proven in Theorem~\ref{thm:symplecticity} (via standard backward-error analysis for symplectic integrators).  The oscillation amplitude scales as ${O(\tau^2\omega_{\max}^2)}$, so methods with larger $\omega_{\max}$ show correspondingly larger oscillations at the same timestep. Importantly, the amplitude reflects the CFL margin, not the quality of the spatial approximation---the eigenvalue analysis of Section~\ref{SPAIOPresults} remains the appropriate metric for operator accuracy.

\subsection{Eigenmode dispersion analysis\label{eigenmodeDispersion}}

A central consequence of different choices of $\mQ_1$ is that the discrete curl-of-curl operator ${\tilde{\mK}=\mQ_1\mP_2}$ has method-dependent eigenvalues, and these eigenvalue errors translate directly into numerical dispersion---i.e., errors in the propagation speed of discrete wave modes. To demonstrate this connection quantitatively, we initialize a time-domain simulation with a single eigenmode of the exact FEEC operator ${\mM_1^{-1}\mP_2}$ and measure the resulting numerical frequency.

The setup follows a common approach of dispersion analysis in FDTD \cite{taflove_computational_2005,blinne_systematic_2018}: propagate a known wave mode and compare its numerical frequency $\hat{\omega}$ to the exact value $\omega$. On the ${20\times20}$ periodic cubical mesh described in Section~\ref{timeDomainEnergy}, we compute the lowest non-kernel eigenmode $\w{v}$ of ${\mM_1^{-1}\mP_2}$ via SLEPc \cite{hernandez_slepc_2005}, obtaining ${\omega=1.0041}$ (corresponding to the fundamental mode of the discrete torus). We initialize ${\w{e}(0)=\w{v}}$, ${\w{b}(0)=\w{0}}$ and evolve for 200 wave periods under Strang splitting with CFL safety factor 0.3. The phase projection ${\langle\w{e}(t),\w{e}(0)\rangle}$ oscillates at the numerical frequency $\hat{\omega}$, which we extract from zero crossings.

For the SPAI-OP run, the probing vector in Eq.~(\ref{SPAIOPobjective}) is set to ${\mP_2\w{v}}$---the curl-of-curl action on the very eigenmode being propagated. With only a single probed mode, the probing weight is set to ${\lambda=5\times10^5}$, which drives the residual error ${\norm{\tilde{\mK}\w{v}-\omega^2\w{v}}/\omega^2}$ on the target mode to below $10^{-4}$ while the overall Frobenius fit remains within $8\%$ of the SPAI baseline---well within the $\alpha=25\%$ tolerance criterion of Section~\ref{lambdaScan}.  (A single-mode probe is easily accommodated by the many degrees of freedom of $\mQ_1$, so the tolerance budget never binds here.)

\begin{table}[h!]
\centering
\caption{Numerical dispersion errors for a single eigenmode on a ${20\times20}$ cubical mesh. The exact eigenfrequency is ${\omega=1.0041}$. SPAI-OP reduces the frequency error to the limit of the zero-crossing extraction method ($<10^{-4}$), achieving essentially exact dispersion accuracy on the probed mode.}
\label{dispersionTable}
\begin{tabular}{lcccc}
\hline
Method & $S(\mQ_1)$ & $\hat{\omega}$ & $(\hat{\omega}-\omega)/\omega$ & PCG iters \\
\hline
Diagonal (Yee) & $S(\mOne)$ & $1.1499$ & $+14.5\%$ & --- \\
SPAI & $S(\mM_1)$ & $0.9634$ & $-4.1\%$ & --- \\
\textbf{SPAI-OP} & $S(\mM_1)$ & $\mathbf{1.00412}$ & $\mathbf{<10^{-4}}$ & $O(10^{2})$ \\
\hline
\end{tabular}
\end{table}

The results (Table~\ref{dispersionTable}) demonstrate that SPAI-OP reduces the frequency error on the targeted eigenmode by more than three orders of magnitude: from $4.1\%$ (SPAI) and $14.5\%$ (diagonal) to below $10^{-4}$---the resolution limit of our zero-crossing fit.  This is a direct time-domain manifestation of the $\sim 10^3\times$ eigenvalue-residual reduction established in Sections~\ref{SPAIOPsubsect} and \ref{SPAIOPresults}: the operator-probed term of Eq.~(\ref{SPAIOPSylvester}) drives $\tilde{\mK}\w{v}\to\omega^2\w{v}$ on the target mode, so the time integrator propagates that mode at essentially its exact frequency. No analogous result is available from standard SPAI or Yee-like lumping, in which the inverse-mass approximation is mode-agnostic.

\begin{figure*}[t!]
\centering
\includegraphics[width=\textwidth]{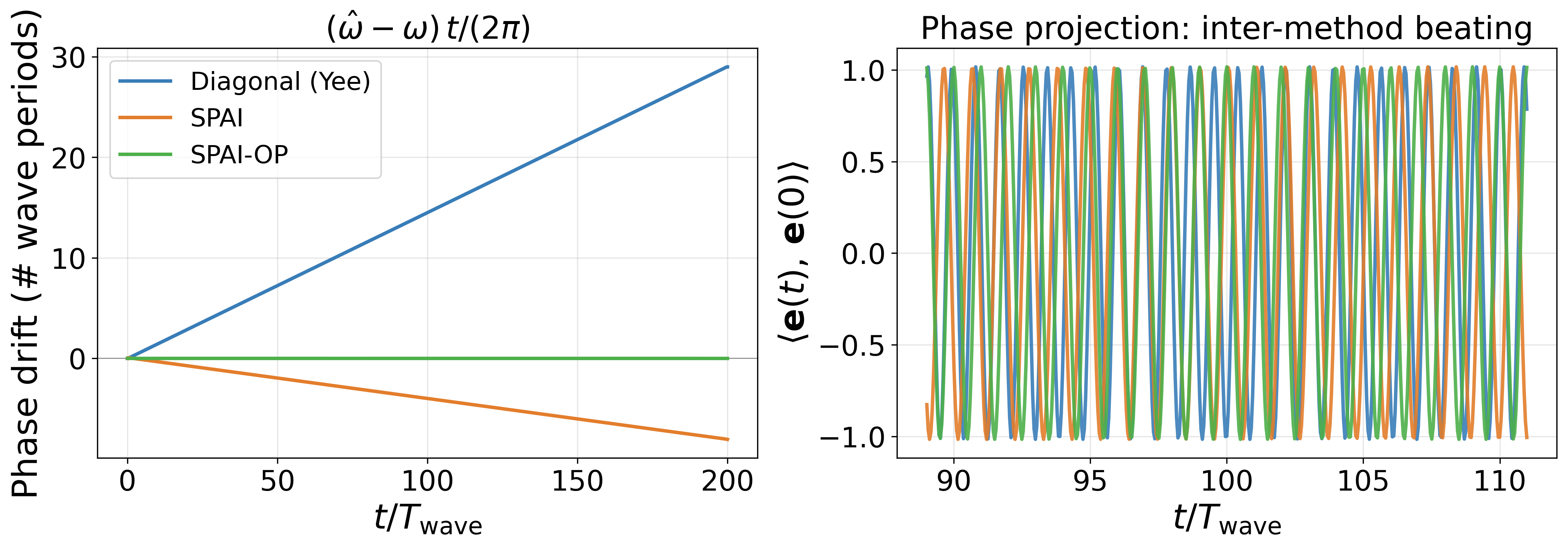}
\caption{Eigenmode dispersion comparison on a ${20\times20}$ cubical mesh over 200 wave periods. \textbf{Left}: Cumulative phase drift ${(\hat\omega-\omega)\,t/(2\pi)}$, in units of whole wave periods. Diagonal (Yee) and SPAI accumulate drift at rates set by their frequency errors of Table~\ref{dispersionTable}; SPAI-OP remains within a small fraction of a wave period over the entire run. \textbf{Right}: Raw phase projection ${\langle\w{e}(t),\w{e}(0)\rangle}$ in a mid-run window. The three methods, initialized in phase, have accumulated distinct numerical phases, producing an inter-method beating pattern.}
\label{dispersionFig}
\end{figure*}

Fig.~\ref{dispersionFig} translates these frequency errors into the time domain. The left panel plots the cumulative phase drift ${(\hat\omega-\omega)\,t/(2\pi)}$, in units of whole wave periods: the diagonal (Yee) signal runs ${\sim}29$ periods ahead of the exact mode by the end of the 200-period run, SPAI falls ${\sim}8$ periods behind, while SPAI-OP remains within a small fraction of a period throughout. The right panel zooms into a mid-run window where the raw phase projection ${\langle\w{e}(t),\w{e}(0)\rangle}$ is plotted directly: the three methods---initialized in phase---have by this point accumulated visibly distinct phases, producing a beating pattern that is a direct time-domain signature of inter-method dispersion error.

The probing weight $\lambda$ controls a continuous tradeoff between overall Frobenius accuracy and the single-mode constraint. At ${\lambda=0}$, SPAI-OP reduces to symmetry-constrained SPAI; as $\lambda$ increases, accuracy concentrates on the probed mode.  At ${\lambda=5\times10^5}$ the probed-mode residual has saturated (any further increase in $\lambda$ buys no visible frequency-error improvement on this mode) while the all-mode Frobenius error is only $8\%$ above the SPAI baseline---comfortably within the $\alpha=25\%$ tolerance budget, since a single-mode probe places only one constraint on the many degrees of freedom of $\mQ_1$.  Practitioners concerned with a narrow band of physical frequencies (typical for dispersion-critical applications such as beam-driven instability simulations) can therefore eliminate dispersion error on the targeted band at essentially one-time setup cost.

\section{Discussion\label{ConclusionSect}}

We have demonstrated that Yee's method can be interpreted as a structure-preserving splitting method within the FEEC formalism, using a cubical mesh with simplified (diagonal) mass matrices. This perspective reveals Yee's algorithm to be a special case of a broader family of methods---generalized Yee methods (GYMs)---summarized by Eq.~(\ref{preYeeFEECEqns}) and Table~\ref{YeeVsGYM}. By Theorem~\ref{thm:symplecticity}, any sparse SPD approximation ${\mQ_1\approx\mM_1^{-1}}$ yields a symplectic and stable GYM, decoupling the choice of mass matrix approximation from the structural guarantees of the method. This decoupling is the key enabling result: it frees the practitioner to optimize $\mQ_1$ for accuracy without sacrificing the conservation properties essential to long-time electromagnetic simulation.

The error scaling analysis of Section~\ref{errorScaling} (Fig.~\ref{CurlOfCurlFig}) establishes that standard Frobenius-optimal SPAI at $\mM_1$-sparsity achieves nearly the full convergence rate of exact FEEC for first-order Whitney forms (${h^{0.73}}$ vs.\ ${h^{0.85}}$), while extending the convergent regime from ${h/\lambda\gtrsim 0.2}$ (diagonal) down to ${h/\lambda\gtrsim .05}$. For second-order elements, ${(\mM_1)^2}$-sparsity captures much of the higher-order improvement (${h^{1.61}}$ vs.\ ${h^{1.98}}$). These results confirm that sparse GYMs can recover substantial finite element accuracy at communication costs comparable to Yee's method.

The SPAI-OP formulation of Section~\ref{SPAIOPsubsect} addresses a fundamental limitation of standard SPAI: its uniform distribution of approximation error across all eigenmodes. By augmenting the Frobenius objective with geometric probing constraints on user-specified wave modes, the resulting symmetry-constrained Sylvester system (Eq.~\ref{SPAIOPSylvester}) sharply concentrates accuracy on those modes. On the 2D structured grid of Section~\ref{SPAIOPresults}, SPAI-OP at the tolerance-optimal $\lambda^*(\alpha=25\%)$ reduces the mean eigenvalue error on the probed modes by factors of $2.2$--$7.8\times$ relative to symmetric SPAI, at a uniform cost of ${\sim}1.3\times$ degradation in the mean error of unprobed modes (Table~\ref{SPAIOPeigTable}). The practitioner can tune this tradeoff continuously through the single parameter $\lambda$ (Fig.~\ref{lambdaScanFig}), from $\lambda=0$ (symmetric SPAI) to the single-mode probe regime where ${\tilde{\mK}\w{v}\to\omega_v^{2}\w{v}}$ on the target mode $\w{v}$ to machine precision. The matrix-free PCG solver with diagonal preconditioning is efficient at the tolerance-optimal $\lambda^*$ values used throughout our tests (${\sim}10^2$ iterations on structured grids), a modest one-time setup cost amortized over the entire time-domain simulation.

The time-domain simulations of Sections~\ref{timeDomainEnergy} and \ref{eigenmodeDispersion} exhibit behavior consistent with symplecticity and demonstrate spectral accuracy in practice. All choices of $\mQ_1$ produce bounded Hamiltonian oscillation over 200 wave periods with no secular drift (Fig.~\ref{energyFig})---a necessary consequence of the symplecticity proven in Theorem~\ref{thm:symplecticity} under backward-error analysis. The eigenmode dispersion test of Section~\ref{eigenmodeDispersion} provides the most direct demonstration of SPAI-OP's value: by probing a specific eigenmode of the curl-of-curl operator, SPAI-OP drives the numerical frequency error below the $10^{-4}$ resolution floor of our zero-crossing diagnostic, relative to $4.1\%$ for symmetric SPAI and $14.5\%$ for diagonal (Yee) (Table~\ref{dispersionTable}, Fig.~\ref{dispersionFig})---essentially eliminating dispersion error on the targeted mode at the cost of a few hundred PCG iterations at setup time. The CFL verification of Section~\ref{CFLverification} independently validates Proposition~\ref{prop:stability}.

Several directions for future work emerge naturally. First, application of SPAI-OP to PIC simulations---via higher-order B-spline bases compatible with smooth particle coupling (Section~\ref{PICsect})---would broaden the applicability of the method. Second, adaptive selection of probing vectors tuned to specific physical regimes (e.g., beam-driven instabilities in PIC, where numerical Cherenkov radiation is particularly problematic \cite{godfrey_cherenkov_2014}) could further enhance the practical utility of SPAI-OP. Third, the symmetry-constrained Sylvester framework itself---as a general technique for computing symmetric sparse approximate inverses of SPD matrices with targeted subspace accuracy, via a scalable matrix-free PCG---may find applications beyond the electromagnetic setting considered here, wherever the accuracy of a sparse inverse must be concentrated in physically important spectral regions.

\section{Acknowledgments\label{AckSect}}

Thank you to Nat Fisch, Ian Ochs, Eli Kolmes, Tal Rubin, Mike Mlodik, Josh Burby and Tyrus Berry for helpful discussions, and to Phil Morrison for his support. This work was supported by the U.S. Department of Energy (DE-AC02-09CH11466), as well as the U.S. Department of Energy Fusion Energy Sciences Postdoctoral Research Program administered by the Oak Ridge Institute for Science and Education (ORISE) for the DOE. ORISE is managed by Oak Ridge Associated Universities (ORAU) under DOE contract number DE-SC0014664. All opinions expressed in this paper are the authors' and do not necessarily reflect the policies and views of DOE, ORAU, or ORISE. This work was further supported by ARPA-E Grant No. DE-AR0001554.

\section*{Credit author statement}
\textbf{Alexander S. Glasser:} Conceptualization, Methodology, Software, Formal analysis, Investigation, Visualization, Writing -- original draft, Writing -- review \& editing.
\textbf{Hong Qin:} Conceptualization, Supervision, Writing -- review \& editing.

\section*{Data availability}
No external datasets were used in this work. The Fortran, Matlab, and Python code generating all figures and tables are available from the corresponding author upon request.

\section*{Declaration of generative AI and AI-assisted technologies in the manuscript preparation process}
During the preparation of this work, the authors used Anthropic's Claude to assist with text editing and development and debugging of numerical code. After using this tool, the authors reviewed and edited the content as needed and take full responsibility for the content of the published article.

\bibliographystyle{IEEEtran}
\bibliography{allrefs.bib}

\end{document}